\documentclass{article}

\usepackage{arxiv}

\usepackage[utf8]{inputenc} 
\usepackage[T1]{fontenc}    
\usepackage{hyperref}       
\usepackage{url}            
\usepackage{booktabs}       
\usepackage{amsfonts}       
\usepackage{nicefrac}       
\usepackage{microtype}      
\usepackage{lipsum}		
\usepackage{graphicx}
\usepackage{natbib}
\usepackage{doi}
\usepackage{hyperref}
\usepackage{algorithm}
\usepackage{algorithmic}

\usepackage{amsmath}
\usepackage{array}
\usepackage{booktabs}
\usepackage{caption}
\usepackage{epstopdf}
\usepackage[figuresright]{rotating}
\usepackage{float}
\usepackage{graphicx} 
\usepackage{longtable}
\usepackage{makecell}
\usepackage{multirow}
\usepackage{natbib}
\usepackage{pdflscape}
\usepackage{ragged2e}
\usepackage{setspace}
\usepackage{subcaption} 
\usepackage{vcell}
\allowdisplaybreaks[4]
\bibpunct[, ]{(}{)}{;}{a}{}{,}

\title{Enhancing efficiency in the face of disruption time: a hybrid memetic-ANS optimization algorithm for the home health care and home care routing and rescheduling problem}

\author{ \href{https://orcid.org/0000-0000-0000-0000}{\includegraphics[scale=0.06]{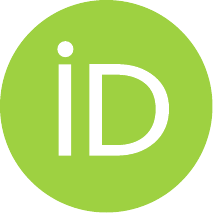}\hspace{1mm}Qiao Pan}\thanks{Use footnote for providing further
		information about author (webpage, alternative
		address)---\emph{not} for acknowledging funding agencies.} \\
	Department of Management and Economics\\
	Tianjin University\\
	Tianjin, 300072 \\
        China\\
	\texttt{panq@tju.edu.cn} \\
	\And
	\href{https://orcid.org/0000-0000-0000-0000}{\includegraphics[scale=0.06]{orcid.pdf}\hspace{1mm}Zhaofang Mao} \\
	Department of Management and Economics\\
	Tianjin University\\
        Tianjin, 300072 \\
        China\\
	\texttt{maozhaofang@tju.edu.cn} \\
}

\renewcommand{\shorttitle}{\textit{Enhancing efficiency in the face of disruption time: a hybrid memetic-ANS optimization algorithm for the home health care and home care routing and rescheduling problem}}

\hypersetup{
pdftitle={Enhancing efficiency in the face of disruption time: a hybrid memetic-ANS optimization algorithm for the home health care and home care routing and rescheduling problem},
pdfsubject={Mathematics, Optimization and control},
pdfauthor={Qiao Pan, Zhaofang Mao},
pdfkeywords={Home health care and home care problem, Rescheduling problem with the rejection of new jobs, Orienteering problem with mandatory visits, Memetic algorithm},
}

\begin{document}
\maketitle

\begin{abstract}
 This paper addresses a realistic home health care and home care (HHC\&HC) problem which has become increasingly complex in the face of demographic aging and post-COVID-19 disruptions. The HHC\&HC sector, as the essential component of modern health care systems, faces unique challenges in efficiently scheduling and routing caregivers to meet the rising demand for home-based care services. Traditional approaches often fall short in addressing the dynamic nature of care requests, especially in accommodating new, same-day service requests without compromising scheduled visits. To tackle these issues, We define the problem as an HHC\&HC routing and rescheduling problem with rejection of new customers (HHC\&HCRRP-RNC), focusing on rescheduling for a single HHC\&HC caregiver in response to new customer requests within a single period. This problem is a variant of both the single-machine reschedule problem and the orienteering problem with mandatory visits (OPMV), where certain nodes must be visited while others are optional. A mixed integer linear programming (MILP) model is developed to cater to two groups of customers: pre-scheduled existing customers and same-day service new customers. The model emphasized maintaining minimal disruptions to the original schedule for existing customers as a constraint, highlighting the balance between adhering to scheduled visits and accommodating new customers. A hybrid memetic-Adaptive Neighborhood Search (ANS) optimization algorithm is proposed to tackle the model. This approach aims to minimize operational costs and opportunity costs while enhancing service quality and patient satisfaction. Through computational experiments, our proposed algorithm demonstrates notable performance, offering significant improvements in both efficiency and robustness within the problem domain.
\end{abstract}


\section{Introduction}\label{introduction}
Health and care provision at home, also commonly referred to as (HHC\&HC), represents home-based services covering health care and personal care based on the degree of medical support required by the customers (also called clients or consumers, hereafter in this paper will be referred to as customers)\citep{markkanenItChangedEverything2021}. HHC\&HC has emerged as a credible substitute for HC\&C in traditional facilities such as hospitals, outpatient clinics, emergency departments, and nursing homes \citep{samaImpactsCOVID19Pandemic2021}. The service range of HHC\&HC can be broad and diverse, typically encompassing nursing care, complex health care procedures, and assistance with daily living activities within their care remit \citep{bazirhaStochasticHomeHealth2023}. Regardless of the variety of services, the ultimate objective is to ensure the prompt delivery of high-quality medical or everyday care from the HHC\&HC center to customers’ homes following their requests \citep{manavizadehUsingMetaheuristicAlgorithm2020a}. Consequently, the primary challenge for HHC\&HC companies is the efficient scheduling of patient appointments and the decisions of optimal routes for HHC\&HC caregivers. These problems are widely recognized in the existing literature as the home health care routing and scheduling problem (HHCRSP) or the home care routing and scheduling problem (HCRSP), which is an advanced variant of the classic vehicle routing problem (VRP), to devise a sequence of routes that optimally balance the financial goals of health care administrators - such as reducing costs, minimizing travel time or distance, and decreasing fuel consumption - and the aim of enhancing patient satisfaction through the provision of more efficient and timely services \citep{AZADEH2015217,cisseProblemsRelatedHome2017}. In this paper, we introduce a composite term - the home health care and home care routing and scheduling problem (HHC\&HCRSP) - to reflect the convergence of these problems from an operations research (OR) standpoint, acknowledging their shared objectives and the operational dilemmas they present in optimizing the delivery of home-based health and care services.

In the aftermath of the COVID-19 pandemic, the health care systems worldwide have been navigating through an era of unprecedented demand, leading to significant strains on resources and amplifying public concerns about infection risks and general hygiene practices \citep{SENCROWE202156,nuraimanDecompositionApproachPrioritised2022}. This phenomenon has underscored the importance of HHC\&HC services, particularly as fears persist about the potential resurgence of the coronavirus and other pandemics. HHC\&HC has become increasingly critical, catering primarily to vulnerable populations such as the elderly, individuals with mild health issues or chronic conditions, and those susceptible to infections. These groups often express a strong preference for receiving health and care services in the safety, cleanliness, and comfort of their own homes, avoiding the potential risks associated with overcrowded traditional health and care facilities \citep{liskerAmbulatoryManagementModerate2021,liu-LargeScalePeriodicHome2020,McKinsey}. 

The HHC\&HC model typically requires booking to facilitate meticulous planning and coordination of resources. Because the services often involve personalized preparation for future treatment or care, which includes medical, paramedical, or social services, it is necessary to arrange early scheduling to match customers with caregivers possessing the appropriate skills and availability \citep{10.1093/fampra/cmac062}. As more customers accustomed to facility-based health and care services now turn to HHC\&HC, traditional approaches to planning and scheduling are increasingly challenged. Decision makers within the HHC\&HC sector are now compelled to make difficult choices in determining which customers to support when integrating the new, often same-day, service requests. This dilemma involves striking a delicate balance between preserving the high quality of service for existing customers and accommodating new requests. The integration of same-day requests requires a careful recalibration of the existing schedules while maintaining high service standards at the expense of frequently declining new customers may discourage them from seeking services in the future, thereby adversely impacting the overall accessibility and utility of HHC\&HC services.  

To tackle these issues, this paper considers rescheduling a single HHC\&HC caregiver to accommodate new customer requests within a single service period, a dilemma we term the HHC\&HCRRP-RNC. Naturally and simultaneously, the problem can be considered as an application in HHC\&HC of the single-machine rescheduling problem with the arrival of new jobs where job rejection is allowed \citep{hallReschedulingNewOrders2004}. Also, it shares similarities with the orienteering problem with mandatory visits (OPMV) which is a novel extension of the Orienteering Problem (OP) as proposed by \citet{palomo-martinezHybridVariableNeighborhood2017}, differentiating certain nodes termed mandatory nodes that must be visited between others called optional nodes that do not require a visit. We develop a MILP model to match the problem considering two distinct customer categories: pre-scheduled existing customers and same-day service new customers (hereafter will be abbreviated as existing customers and new customers, respectively). Assuming optimal initial scheduling for existing customers to minimize total travel costs, these customers are guaranteed specific visit times with fairly minor deviations following the schedule. Hereby, we formulate an original problem where the objective is to determine the most efficient route for a single HHC\&HC caregiver while adhering to the promised timeframes for existing customers. Accordingly, the original problem can be seen as a Traveling Salesman Problem with soft Time Windows (TSPsTW). 

When introducing new customers who may arrive unexpectedly and potentially disrupt the pre-established schedule, it necessitates a decision-making process on whether to accept or reject these requests based on the caregiver's capacity and the potential impact on scheduled visits. This process is influenced by two key performance indicators: firstly, the total income from accepted customers; and secondly, the costs associated with rejecting new requests. The latter represents a dual loss, encompassing both the immediate financial loss and the long-term opportunity cost of permanently losing potential customers. Furthermore, our model treats the delay imposed on existing customers, namely disruption time, as a constraint rather than an objective component to emphasize the importance of adhering to the original schedule for existing customers.

Given that the complexity of rescheduling and the orienteering problem (OP) are both identified as NP-hard \citep{hallReschedulingNewOrders2004,https://doi.org/10.1002/1520-6750(198706)34:3<307::AID-NAV3220340302>3.0.CO;2-D}, the HHC\&HCRRP-RNC is similarly classified as NP-hard. Several heuristics and exact algorithms in the literature have been successfully addressed rescheduling and OPMV. Evolutionary algorithms (EAs) have demonstrated effectiveness in a broad spectrum of scheduling problems both in single and multi-objective machine scheduling contexts \citep{BALLESTIN20082175}. Their success may be attributed to the ability to handle complex, multi-faceted problems through adaptive and iterative processes. In terms of OPMV, to the best of our knowledge, only one exact has been proposed, i.e. a branch-and-check algorithm of \citet{vuBranchandcheckApproachesTourist2022}. Heuristics applied in OPMV consist of adaptive large neighborhood search (ALNS)\citep{angelicasalazar-aguilarMultidistrictTeamOrienteering2014,fangRoutingUAVsLandslides2023}, simulated annealing (SA)\citep{linSolvingTeamOrienteering2017}, variable neighborhood search (VNS)\citep{palomo-martinezHybridVariableNeighborhood2017}, neural network (NN)\citep{fangRoutingUAVsLandslides2023}, and notably, the memetic algorithm (MA)\citep{luMemeticAlgorithmOrienteering2018}, which has been identified as particularly potent in this context. The MA's strength lies in its integration of tabu search (TS) and genetic algorithms (GAs), showcasing superior performance through a blend of global and local search strategies. 

In addressing the HHC\&HCRRP-RNC, this paper advances the application of MA by \citet{luMemeticAlgorithmOrienteering2018}, employing a different framework by incorporating ALNS in the initial solution phase and introducing unique move operators within the TS local search phase. This adaption is designed to meet the specific demands of the HHC\&HCRRP-RNC, with ALNS enhancing the diversity of the initial solution population and TS refining solutions through targeted local searches. Note that although ALNS is not typically classified as a population-based algorithm, this innovative combination aligns with the MA's principle of integrating various search methodologies to optimize problem-solving processes \citep{MemeticAlgorithms}, setting a new benchmark for addressing the complex rescheduling challenges inherent in HHC\&HC domains. 

The remainder of this paper is organized as follows. In Section \ref{sec_literature}, we briefly review the relevant studies. The HHC\&HCRRP-RNC is formally described in Section \ref{sec_model}, along with a MILP formulation and the notation adopted throughout the paper. Our proposed MA is presented in Section \ref{sec_algorithm}. In Section \ref{sec_experiment}, the results of the MA are presented, including statistical comparisons with existing state-of-the-art exact and heuristic methods. Finally, the paper concludes with Section \ref{sec_conclusion}, where we summarize our findings and implications of the research. 

\section{Literature review}\label{sec_literature}

In this section, we review the relevant and recent studies in two streams: (i) the rescheduling and scheduling problem with rejection allowed; and (ii) routing and scheduling for HHC\&HC.

\subsection{Research on rescheduling and scheduling with rejection allowed}
In typical scheduling scenarios represented by pull-based manufacturing line optimization, disruptions are often unavoidable due to constraints involving downstream demand-side order and technical standard requirements, production-side machine status and staffing, and upstream supply-side raw material supplies. When disruptions occur, it is essential to arrange a new scheduling plan based on current production needs and production factors. \citet{hallReschedulingNewOrders2004} first introduced the concept of rescheduling, constructing models based on a single-machine environment with the total number of tardy jobs and the total tardiness of jobs as objectives, and setting upper limits for these two parameters in the constraints. However, this study only considered the total number and time of tardy jobs when setting the objective function and did not take into account other classic scheduling costs. Moreover, the problem background of the study was relatively singular and did not involve common scenarios in the scheduling field such as flow shops and job shops. \citet{wang2017multi} constructed scenarios where the due time changes temporarily and jobs can be rejected, proposing a rescheduling problem with objectives of minimizing completion time, keeping the perturbation time within limits, and not exceeding the total rejection cost. At the same time, this study used dynamic programming to verify that when the number of machines is fixed, the problem is NP-hard, thereby demonstrating the significance of developing heuristic algorithms for solving medium to large-scale rescheduling problems. \citet{luo2018rescheduling} studied two types of rescheduling problems: the problem of minimizing the maximum weighted completion time with a maximum disruption time constraint and relaxing the maximum disruption time constraint. \citet{rener2022single} solved the rescheduling problem aimed at minimizing the total absolute disruption time using branch and memorize (an exact algorithm), but the initial upper bound in the study was still generated using a heuristic algorithm. \citet{zhang2022rescheduling} constructed a rescheduling problem with the objective function of minimizing weighted tardiness and disruption time as a constraint. The study set two types of disruption time constraints: limiting the total disruption time and limiting the disruption time per occurrence; the study proved that the minimization of weighted tardiness under both constraints is an NP-hard problem, hence heuristic algorithms are more suitable for solving medium to large-scale examples.

\subsection{Research on routing and scheduling for HHC\&HC}
The numerous existing studies on HHC\&HC have proposed a series of models based on different problem formulations, reflecting a wide range of features associated with different real-world scheduling scenarios. Yet most of them focus on the broader challenges of scheduling in HHC\&HC, and there is a paucity of contributions related to the complexities arising from dynamic events \citep{fikar2017home}. In reality, this field often involves highly uncertain and diverse elements, such as the travel time of medical staff to different patients, the duration of care, and the provision of medical services being variable. Different medical personnel often have different skills and specialties, and patients' preferences for the skills, gender of the medical staff, and service time vary widely \citep{gutierrez2013home}. Therefore, how to efficiently allocate medical personnel (both staff and time), design scientific visit routes to maximize the conservation of medical human resources, and save costs while respecting patient needs, has become a focus in the field of HHC\&HC services research. In a typical HHCRSP, before the service starts, the center arranges the visiting schedule of medical staff based on the patient's location and personal needs. HHC\&HC staff usually depart from the center, use vehicles provided by the center to go to the patients' homes that have been arranged in advance, and return to the center at the end of the day's work. Since HHCRSP and classic VRP have many common points in terms of objectives and important problem characteristics (for example, both aim to (i) arrange service vehicles for customers located in different places; (ii) minimize the vehicle travel time or distance; important characteristics include: (i) each customer must be visited once and only once by one of the vehicles; (ii) the start and end points of the route are both the center; etc.), HHCRSP has been classified as a variant of VRP in many previous studies. Additionally, HHCRSP shares many similarities with many common VRP extension problems, such as the customer time window constraints of VRPTW are also common in the settings of HHCRSP \citep{di2021routing}.

Other models proposed in the literature to handle such uncertainties in HHC\&HC are from studies of \citet{koeleman2012optimal}, \citet{carello2014cardinality}, \citet{rodriguez2015staff} and \citet{martinez2019re}. The Markov decision model \citep{koeleman2012optimal} is introduced to address personnel planning in care-at-home service facilities within a stochastic framework. This model is discussed in contexts both with and without waiting rooms, as well as proves some monotonicity properties of the value function and the structure of the optimal policy, while the complexity of the model makes it less accessible and user-friendly for decision-makers in real-world settings. A two-stage stochastic programming model \citep{rodriguez2015staff} is employed to determine the minimal human resource requirement that can adequately cover a specified percentage of operational days. This research primarily concentrates on the uncertainties in demand stemming from various pathologies and locations, thereby characterizing the research within the realm of staff dimensioning problems. A decomposition approach \citep{martinez2019re} is proposed to reduce HHC worker idle time while adhering to continuity of care constraints. The approach is specifically designed for situations where an existing patient leaves the HHC agency, as opposed to scenarios where a new patient is about to enter. In the robust optimization model \citep{carello2014cardinality}, a different scenario is contemplated where each patient consistently receives service from the same nurse, creating a nurse-patient continuity problem. However, there are computational challenges associated with this approach. Also, the model allows patients’ demands to take only two values, i.e. expected and maximum. This restriction might result in overly conservative solutions.

\section{Problem description and modeling} \label{sec_model}
In this section, we provide a detailed formal description and a MILP model for the HHC\&HCRRP-RNC. Initially, we introduce the original problem that only existing customers are considered, leading to the formulation of an original schedule. Subsequently, we present the rescheduling model which incorporates the arrival of new customer requests, using the original schedule as the baseline for adjustments. 

\subsection{Problem description and modeling}
The HHC\&HCRRP-RNC is defined on a complete undirected graph $G=(V, A)$, where each node (i.e., customers' home locations) $i\in V$ is associated with a non-negative score $p_i$ representing the payment received from customers who are serviced. The node set \textit{V} is divided into $V^E=\{0\}\cup V^E$ and $V^N=\{0\}\cup V^N$, which present the existing and new customers' home locations, respectively. $A=\{(i,j)|i,j\in V,i\neq j\}$ is the arc set that connects the pairs of nodes, and the depot is defined as node 0. For simplicity in calculations without compromising the essence of the problem, the required time to travel from a node \textit{i} to \textit{j}, which is denoted as $t_{ij}$, is assumed to be directly equivalent to the Euclidean distance between the two nodes and the triangle inequality holds. The requests of the existing customers are recorded before the planning horizon, allowing the decision-maker to generate an original schedule $\pi^*$, which serves as a baseline and sets a foundation for subsequent rerouting and rescheduling decisions within the predefined maximal disruption time \textit{T}. When new customers are added, the deviation in arrival times at each existing patient \textit{i} between the original $\pi^*$ and rescheduled plan $\pi$ should not exceed \textit{T}, i.e. $\Delta_i(\pi)=|s_i(\pi)-s_i(\pi^*)|,i\in V^E$, ensuring minimal disruption to the original schedule. Additionally, if any new customers are rejected in the rescheduling model, a penalty \textit{r} is incurred for each rejection case. This penalty accounts for not only the immediate loss of payment but also the opportunity cost associated with the potential permanent loss of these customers. The inclusion of this penalty in the model underscores the economic and relational impact of customer rejection in the rescheduling process.

\subsection{Assumptions}
In summary, based on the previous description, the following assumptions are considered for our proposed models:

\begin{enumerate}
	\item Without losing any generalization, the planning horizon has not been constrained to a specific time duration.
	\item For reduction of complexity, the proposed model in this paper is a slightly simplified version of HHC\&HCRSP by omitting the service time. Thus, the arrival and departure time at each served patient's home location should be the same. This simplification facilitates a more focused analysis of other critical dimensions.
	\item The HHC\&HC caregiver starts the route at the HHC\&HC organization (it served as the depot, hereafter will be abbreviated as HHC\&HCO),  serving all customers in the rescheduled plan and returns to the HHC\&HCO.
	\item Each customer can be served at most once.
	\item Each existing customer must be served exactly once.
	\item The caregiver is allowed to reject any customers, if necessary. Each rejection case will incur a fixed penalty.
	\item Deviation of arrival time at each existing customer's home location should not exceed the maximal allowed disruption time between the original and rescheduled plan.
	\item Transportation costs are directly equivalent to the total Euclidean distances traversed. 
\end{enumerate}

\subsection{Notation}
Based on the above assumptions and before presenting the mathematical formulation of the HHC\&HCRRP-RNC, we outline notations in Table \ref{tab1}.

\begin{table}[htpb]
    \centering
    \setlength{\abovecaptionskip}{0pt}
    \setlength{\belowcaptionskip}{5pt}
    \caption{\normalsize{Notation for the HHC\&HCRRP-RNC formulation}}
    \begin{tabular}{p{5cm}<{\raggedright\arraybackslash}p{10cm}<{\raggedright\arraybackslash}}
        \specialrule{0.05em}{3pt}{3pt}
        Notation & Definition\\ 
        \specialrule{0.05em}{3pt}{3pt}
        \multicolumn{2}{l}{\textbf{Sets:}}\\
        \textit{V} & Set of all customer nodes ($i,j=0,v_1,\cdots,v_{n_E},v_{n_E+1},\cdots,v_{n_E+n_N}$)\\
        $V^O=\{v_1,\cdots,v_{n_O}\}$ & Set of the existing customer nodes ($i,j=v_1,\cdots,v_{n_E}$)\\
        $V^N$ & Set of the new customer nodes ($v_{n_E+1},\cdots,v_{n_E+n_N}$)\\
        \textit{A} & Set of arcs connecting nodes\\
        \multicolumn{2}{l}{\textbf{Parameters:}}\\
        $\pi^*$ & An optimal permutation of nodes for the existing customers\\
        $t_{ij}$ & Distance covered to reach from a customer node \textit{i} to \textit{j}\\
        $\Delta_i(\pi)$ & The disruption time of an existing customer \textit{i} between the original and rescheduled plans $\pi$ and $\pi^*$, $\Delta_i(\pi)=|s_i(\pi)-s_i(\pi^*)|\quad i\in V^E$\\
        \textit{T} & The threshold of disruption time, i.e., $\Delta=max\{\Delta_i\}\le T$\\
        $T_{max}$ & The maximum total traveling time in the reschedule plan\\
        $p_i$ & Payment of a customer \textit{i} ($i\in V^O\cup V^N$)\\
        \textit{r} & Cost of rejecting a new customer\\
        $s_i(\pi^*)\quad i\in V^O$ & The arrival time of the caregiver at an existing customer node $i\in V^E$ in $\pi^*$; hereafter this parameter will be abbreviated as $s_i^O$\\
        \textit{M} & A big number\\
        \multicolumn{2}{l}{\textbf{Decision and auxiliary variables}}\\
        $x_{ij}$ & A binary variable that equals to 1 if the caregiver services arc $(i,j)\in A$, otherwise it equals to 0\\
        $y_i$ & A binary variable that equals 1 if customer $i\in V^O\cup V^N$ is serviced; otherwise it equals to 0\\
        $v_i$ & A binary variable that equals 1 if a new customer $i\in V^N$ is serviced; it equals 0 if that customer is rejected\\
        $s_i$ & The arrival time of the caregiver at a customer node $i\in V^O\cup V^N$ in the rescheduled plan\\
        $u_i$ & The position of node $i\in V$ in the route, in particular, we have $u_0 = 0$ and $u_{n_E+n_N+1}=n_E+n_N+1$\\
        \specialrule{0.05em}{2pt}{0pt}
    \end{tabular}
\label{tab1}
\end{table}

\subsection{The first model: the original problem} \label{sec_original}
The proposed MILP for the original problem is given as follows, which considers only existing customers:

\begin{equation}
    \text{Minimize}\; \sum_{(i,j)\in A}t_{ij}x_{ij}
\end{equation}
\begin{gather}
    \sum_{i\in V}x_{ij}=1\qquad\forall j\in V^E\\
    \sum_{j\in V }x_{ij}=1\qquad\forall j\in V^E\\
    \sum_{i\in V^E}x_{0i}=1\\
    \sum_{i\in V^E}x_{i0}=1\\
    s_i+t_{ij}x_{ij}\le s_j+M(1-x_{ij})\qquad\forall i,j\in V^E,i\neq j\\
    u_i-u_j+1\le n_E(1-x_{i,j})\qquad\forall\ i,j\in V^E,i\neq j\\
    1\le u_i\le n_E-1\quad\forall i=1,\cdots,n_E-1\\
    x_{ij}\in\{0,1\}\qquad\forall i,j\in\{0\}\cup V^E,i\neq j\\
    s_i\ge0\qquad\forall i\in\{0\}\cup V^E
\end{gather}

As we can see, the above formulation is a TSP model, except that we do not impose a maximal total travel time. The objective function (1) is to minimize the sum of the total travel costs (interpreted as time). Constraints (2) and (3) function dually as network constraints and as a restriction to ensure that each customer receives service precisely once. Constraints (4) and (5) ensure that the caregiver commences and ends the route at the HHC\&HCO, which acts as the depot. Constraint (6) forces arrival times at any given nodes must precede those at their immediate successors. Constraints (8) and (9) present the Miller-Tucker-Zemlin formulation to allow the elimination of sub-tours. Finally, constraints (9) and (10) are related to the decision variables.        
\subsection{The second model: the rescheduling problem} \label{sec_rescheduling}
Based on the original schedule as the baseline for adjustments, we present the rescheduling model which incorporates the arrival of new customer requests. The proposed MILP for the rescheduling problem is given as follows:

\begin{equation}
    \text{Minimize}\; \sum_{i\in V^E\cup V^N}p_iy_i-r\sum_{j\in V^N}v_j
\end{equation}
\begin{gather}
    \sum_{i\in V^E\cup V^N}x_{0i}=1\\
    \sum_{i\in V^E\cup V^N}x_{i0}=1\\
    y_i=1\qquad\forall i\in V^E\\
    \sum_{i\in V}x_{ij}=y_j\qquad\forall j\in V^N\\
    \sum_{j\in V}x_{ij}=y_i\qquad\forall i\in V^N\\
    \sum_{i\in V}x_{ij}+v_j=1\qquad j\in V^N\\
    s_i-s_i^O\le z\qquad i\in V^E\\                             
    s_i^O-s_i\le z\qquad i\in V^E\\
    z\le T\\
    \sum_{(i,j)\in A}t_{ij}x_{ij}\le T_{max}\\
    s_i+t_{ij}x_{ij}\le s_j+M(1-x_{ij})\qquad\forall(i,j)\in A\\
    u_i-u_j+1\le(n_E+n_N)(1-x_{ij})\qquad\forall i,j\in V^E\cup V^N\\
    1\le u_i\le n_E+n_N\quad\forall i=1,\cdots,n_E+n_N\\
    x_{ij}\in\{0,1\}\qquad\forall(i,j)\in A\\
    y_i\in\{0,1\}\qquad\forall i\in V^E\cup V^N\\
    v_i\in\{0,1\}\qquad\forall i\in V^N\\
    s_i\ge0\qquad\forall i\in V
\end{gather}

In the second model, objective (11) maximizes the total net profits, i.e. the total payments received from the customers who are served less the costs incurred from the rejection of new customers. Constraints (12) and (13) force the route to start and end at the HHC\&HCO (the depot). Constraint (14) guarantees that each existing customer is served exactly once. Constraints (15) and (16) are network constraints. Constraint (17) creates a binary condition where each new customer is either rejected or served for one time. Constraints (18)-(20) impose the maximal allowed disruption time \textit{T} to arrival time at each existing customer's home location, compared to that in the original plan. Constraint (21) ensures that the total traveling time does not exceed the time budget $T_{max}$. Constraint (22) keeps track of the timeline at each node. Constraints (23) and (24) incorporate the Miller-Tucker-Zemlin formulation to eliminate sub-tours within the routing solution. The remaining constraints are related to definitions of variables.

As shown, the above formulation is a special case of the OPMV model, which was proposed by \citet{luMemeticAlgorithmOrienteering2018} built upon \citet{palomo-martinezHybridVariableNeighborhood2017}. The notable exceptions are: first, our objective function incorporates penalties to quantify the immediate and potential losses incurred from the rejection of new customers. Then, we have integrated constraints that account for the disruption in arrival times at the locations of existing customers. Lastly, unlike the original OPMV model, our formulation does not necessitate the set of incompatible nodes.

\section{Solution method for the HHC\&HCRRP-RNC} \label{sec_algorithm}
In this section, we develop a MA algorithm developed for the HHC\&HCRRP-RNC. The MA is an advanced metaheuristic, a hybrid optimization method that combines elements of both global and local search strategies. The term "memetic" in MA, named by \citet{moscato1989evolution}, draws from the concept of a "meme", analogous to a gene in biology but within the context of cultural information transfer \citep{dawkins1976hierarchical}, metaphorizing the algorithm's ability to evolve solutions not just through genetic mechanism but also via learning and adaption. The MA gathers both the strengths of population-based and local search approaches to enhance a promising dual capability for exploitation and exploration. This equips the MA to navigate the equilibrium between intensive search and the preservation of diversification \citep{zhang2015memetic}. Compared to previous studies, our MA replaces the traditional genetic architecture and the best-position-insertion operator with ALNS. This adjustment allows ALNS to function dually as the initial solution generator and the global search mechanism. It is executed based on our findings, as detailed in Section \ref{sec_experiment}, indicating that ALNS outperforms GA and the best-position-insertion operator in the context of the HHC\&HCRRP-RNC. In addition, we restructure the ALNS framework and employ innovative, customized operators specifically designed to accommodate the tightly constrained problem that differs from traditional OPMV. Also, specially tailored components in TS are integrated. The general framework of our MA for HHC\&HCRRP-RNC is outlined in Algorithm \ref{alg1}.

\begin{algorithm}[h]
	\caption{Memetic algorithm for the HHC\&HC-RNC}
	\label{alg1}
	\begin{algorithmic}[1]
		\REQUIRE \emph{$G=(V,E)$}: a weighted undirected graph
  		\ENSURE $S^*$: best solution found so far
		\STATE $S_0\leftarrow$ \emph{Initial\_Solution} \quad/*the best path of the original problem*/
		\FOR{$i=\{1,2,\cdots,p\}$} 
		\STATE $P=\{S_1,S_2,\cdots,S_p\}\leftarrow$ \emph{ALNS($S_0$)}\quad/*Section 4.1*/
        \ENDFOR
        \STATE $S^*\leftarrow$ \emph{Best($P$)}\quad/*record the best solution found so far*/
        \WHILE{stop condition is not met}
        \FOR{each $S$ from $P$}
        \STATE $S_0 \leftarrow$ \emph{Tabu\_Search($S$}) \quad/*Section 4.2*/
        \IF{$f(S_0)>f(S^*)$}
        \STATE $S^*\leftarrow S_0$
        \ENDIF
        \STATE $P\leftarrow$ \emph{Population\_Update($P,S_0$)}
        \ENDFOR
        \ENDWHILE
	\end{algorithmic}
\end{algorithm}

\subsection{ALNS: finding the initial solution}
\label{ALNS}
In this section, we present an ALNS procedure leveraging three groups of move operators - insertion, internal, and removal- to generate neighborhood solutions. This ALNS procedure is crucial for thoroughly exploring the search space to forge a superior initial population. Traditionally, ALNS sequentially employs a pair of operator groups, typically known as destroy ($\Omega^-$) and repair ($\Omega^+$), to iteratively modify feasible solutions. Each iteration commences with the removal of certain nodes (destroy), followed by a restructuring of the solution (repair), utilizing nodes from the removal list. The algorithm dynamically adjusts the weights of these operators based on their efficacy, selecting the more successful operators for subsequent iterations. However, our problem context is an OP in nature, which significantly differs from TSP as it allows the rejection of new customer nodes. This capability results in feasible solutions of varying lengths after each iteration. In essence, the number of nodes on each output path might differ following each move. Therefore, the traditional logic of pairing destroy with repair moves does not apply to our problem. In addition, given that the initial feasible solution for rescheduling is derived from the original problem's optimal solution - which excludes any new customer nodes -  it is natural to perform a node insertion move in the first step. This will encourage to generation of practical neighborhoods in the following iterations. 

\begin{figure}[h]
	\centering
	\includegraphics[width=16.5cm]{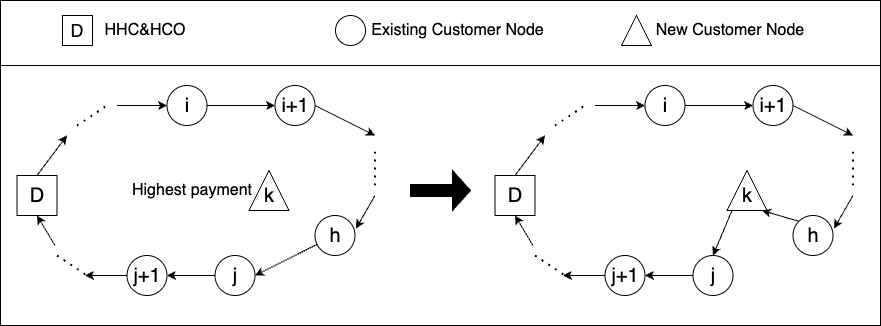}
	\caption{An example of the Restricted Shortest Insertion operator, in which the unvisited node k with the highest payment is inserted into a position that increases the minimal increment of travel time}
	\label{fig1}
\end{figure}

\begin{figure}[h]
	\centering
	\includegraphics[width=16.5cm]{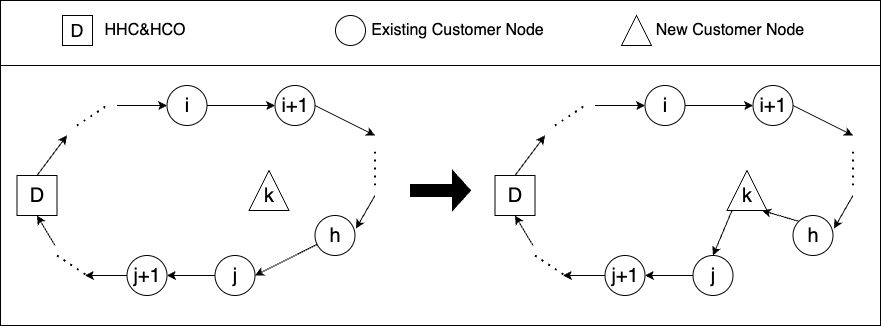}
	\caption{An example of the Restricted Shortest Insertion operator, in which the unvisited node k is inserted so that the least increase in travel time is produced}
	\label{fig2}
\end{figure}

The above description explains why we redesigned and renamed the classical ALNS's destroy and repair operators as new customer node removal and insertion operators, still denoted by $\Omega^-$ and $\Omega^+$ respectively. Also drawing upon the aforementioned observations and as demonstrated within the objective function (see Section \ref{sec_rescheduling}), we prioritize the insertion of as many new customer nodes as possible to concurrently increase the total payments and mitigate penalties associated with new customer rejection. Hence, we provide three new customer node insertion operators as follows:

\begin{enumerate}
    \item Payment-prioritized Insertion: this operator aims to insert the unvisited node \textit{k}k with the highest payment value into an existing arc $(h,j)$ within the current route \textit{S} such that the increment of travel time is minimized (see Fig. \ref{fig1}). To prevent from staying in a local optimum, this operator allows for the insertion of a random unvisited node at a position that results in the minimal increment of traveling time, with a specified probability \textit{P}.
    \item Restricted Shortest Insertion: this operator inserts an unvisited node $k$ into the current route \textit{S} by breaking an arc $(h,j)$ that leads to the least increase in travel time \citep{yu2019matheuristic}. This operator also permits the occasionally insert insertion of a random unvisited node with a probability \textit{P}, at a location within the current route that minimized the increment in traveling time (see Fig. \ref{fig2}).
    \item Disturbed Shortest Insertion: this operator functions similarly to Restricted Shortest Insertion, except that a random value \textit{U} is added to the shortest length of the current route to create a perturbation. Then, with probability $1-P$, it chooses the new shortest arc after perturbation to insert the random new patient node; and with probability \textit{P}, randomly inserts a new patient node on any arc of the current route. More precisely, the perturbation value \textit{U} is computed as $U=L(h,j)\times u\times r$, where $L(h,j)$ represents the length of the arc $(h,j)$ which possesses the shortest length in the current route;  \textit{u} is the perturbation factor and \textit{r} is a uniform $(-1,1)$ random number. 
\end{enumerate}

\begin{figure}[h]
	\centering
	\includegraphics[width=16.5cm]{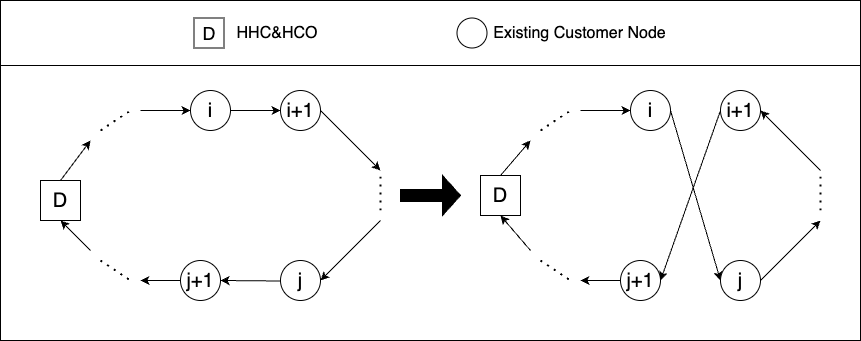}
	\caption{An example of the 2-opt operator which replaces two arcs $(i,i+1)$ and $(j,j+1)$ with $(i,j)$ and $(i+1,j+1)$}
	\label{fig3}
\end{figure}

\begin{figure}[h]
	\centering
	\includegraphics[width=16.5cm]{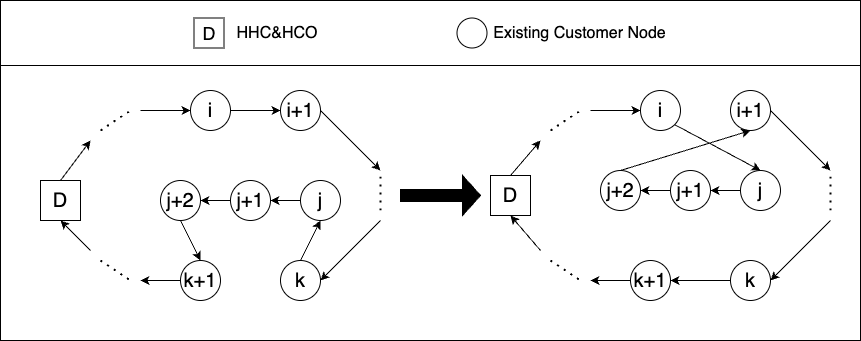}
 	\caption{An example of the Restricted Longest Removal operator in which three consecutive nodes $j$, $j+1$, and $j=2$ are relocated to a different position while keeping the sequence within this chain}
	\label{fig4}
\end{figure}

\begin{figure}[h]
	\centering
	\includegraphics[width=16.5cm]{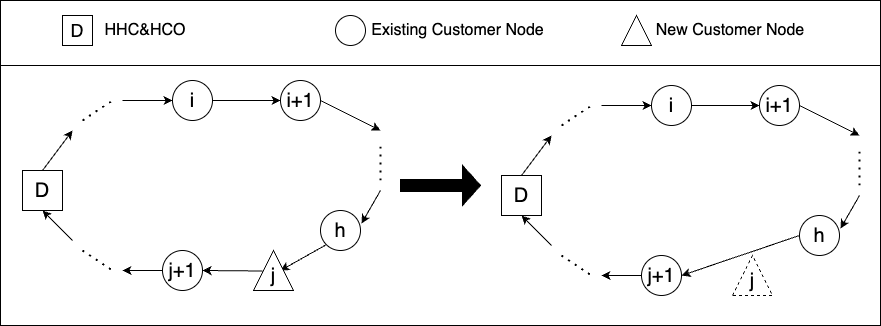}
	\caption{An example of the Restricted Longest Removal operator which removes the visited new customer node $j$ from the current route to produce the most decrease in travel time}
	\label{fig5}
\end{figure}

In the subsequent stage, as in \citet{yu2019matheuristic}, we also employ the 2-opt and Or-opt operators, which are categorized as internal move operators. These two operators are applied to shorten the length of the current route. Additionally, the remaining two operators are classified as new customer node removal operators, functioning specifically with the execution of removing new customer nodes from the current route. The details of internal move and new customer node operators are defined as follows:

\begin{enumerate}
    \item 2-opt: this operator is part of our internal move operators, which are employed solely to reorder the existing nodes in the current route without adding or removing nodes. It replaces two arcs $(i,i+1)$ and $(j,j+1)$ with $(i,j)$ and $(i+1,j+1)$, effectively reversing the sequence of the nodes between $i+1$ and $j$ int the current route (see Fig. \ref{fig3}).
    \item Or-opt: this operator entails relocating a sequence of one, two, or at most three consecutive nodes to a different position within the current route while keeping the original order of this chain (see Fig. \ref{fig4}). 
    \item Payment-prioritized \& Distance-normalized Removal: this operator removes a new customer node \textit{k} from the current route \textit{S}. The condition of removing a new customer node is based on the ratio of payment received from that node to the incremental travel time caused by including that node in the route. Specifically, a node \textit{k} is eligible for removal if and only if the ratio of its payment to its associated increment in travel time is no more than the minimum ratio of payment to increment of travel time for all new customer nodes such that $\frac{Payment_k}{\Delta TravelTime_k}\le\min_{i\in V^N}\frac{Payment_i}{\Delta TravelTime_i}$, where $\Delta TravelTime_i$ is the increment of the travel time between two routes before and after the inclusion of node \textit{i}. In this way, we anticipate eliminating a new customer that contributes less payment while generating more travel distance. 
    \item Restricted Longest Removal: this operator simply removes any new customer node that results in the maximum increment of travel time (see Fig. \ref{fig5}). 
\end{enumerate}

Our ALNS procedure can be outlined in Algorithm \ref{alg2}.

\begin{algorithm}[h]
	\caption{ALNS with insertion, internal, and removal operators}
	\label{alg2}
	\begin{algorithmic}[1]
		\REQUIRE $S_0$: the initial\_Solution which is the best path of the original problem; $\Omega^+$: the set of insertion operators; $\Omega$: the set of internal move operators; $\Omega^-$: the set of removal operators; $L_{ALNS}$: search depth; $u_w$: frequency of updating weights of operators in $\Omega^+$, $\Omega$ and $\Omega^-$; $s^b$: the score of finding a solution surpasses the best solution found so far; $s^c$: the score of finding a solution surpasses the current solution; $s_w$: the score when accepting a worse solution; $P$: the probability to accept a worse solution; $\rho$: the evaporation factor to prevent the solution from converging too quickly on a sub-optimal path; $u$: a perturbation factor to explore more diversified neighborhoods 
  		\ENSURE $S^*$: best solution found so far
		\STATE \emph{Iter(ALNS)} $\leftarrow0$ \quad/*iteration counter*/
        \STATE $S^*\leftarrow S_0$ \quad/*feasible solution with the best objective value*/
        \STATE $S^c\leftarrow S_0$ \quad/*current solution*/
        \WHILE{\emph{Iter} < $L_{ALNS}$}
        \STATE Adaptively select insertion, internal and removal operators $\omega^+$, $\omega$ and $\omega^-$ from $\Omega^+$, $\Omega$ and $\Omega^-$
        \STATE $S \leftarrow \omega^+,\omega,\omega^-(S_0)$\quad/*a feasible solution*/
        \IF{$f(S)>f(S^c)$}
        \STATE $S^c\leftarrow S$
        \IF{$f(S)>f(S^*)$}
        \STATE $S^*\leftarrow S$
        \STATE $Score(\omega^+, \omega ,\omega^-)\;+=\;s_b$
        \ELSE
        \STATE $Score(\omega^+, \omega ,\omega^-)\;+=\;s_c$
        \ENDIF
        \ELSIF{accept $S$}
        \STATE $Score(\omega^+, \omega ,\omega^-)\;+=\;s_w$
        \ENDIF
        \IF{\emph{Iter} mod $u_w=0$}
        \STATE Update Weight($\Omega^+,\Omega,\Omega^-$)
        \STATE Score($\Omega^+,\Omega,\Omega^-)\leftarrow0$
        \STATE Times($\Omega^+,\Omega,\Omega^-)\leftarrow0$
        \ENDIF
        \ENDWHILE
	\end{algorithmic}
\end{algorithm}

The ALNS procedure starts with an initial solution $S_0$, which is the best route found for the original problem (lines 2 and 3). Next, we repeat a loop of the operations between lines 5 and 21 until a stop criterion is met. At each iteration, a new solution $S$ is created (line 9) using insertion, internal, and removal operators sequentially (line 6). It should be noted that each $\omega^+,\omega$, and $\omega^-$  are selected using a weight self-adaptive mechanism, which is based on its effectiveness (line 5). Operators that yield better results (higher scores) relative to their frequency of use are given higher weights, suggesting a dynamical optimization strategy that prioritizes efficient and effective operators. From lines 7 to 16, according to the quality of $S$, if $S$ is feasible, the scores of any operator $\omega^+,\omega$, and $\omega^-$ involved are updated. Initially set to 0, the score is incremented by $s_b$ if $S$ is the new best solution, or by $s_c$ if it is only better than the previous one. A worse new solution can be accepted in our algorithmic settings at a probability $P$ (line 15), and the score is incremented by $s_w$ in this case (line 16). In lines 18 and 19, all operators' weights are updated at every $u_w$ iteration using their accumulated scores and frequencies of usage. After that, the scores and usage frequencies of the operators that are accumulated over the previous iterations, are reset to 0 for the subsequent cycles of the ALNS procedure, as illustrated in lines 20 and 21. The use of such initialization prevents the procedure from abusing the prioritized operators. 

Additionally, an operator's weight is calculated according to Equation \ref{equ29}: 

\begin{equation}
	W_i=\left\{\begin{array}{lr}
		\displaystyle\emph{$W_i\times \rho$}, \quad \emph{$t_i=0$} \vspace{1.5ex}\\
		\displaystyle\emph{$W_i\times(1-\rho)$}+\frac{\emph{$s_i\times\rho$}}{\emph{$t_i$}}, \quad  \emph{$t_i$} \neq 0
	\end{array}\right. \label{equ29}
\end{equation}

where $W_i$ is the current weight of operator $i$, $t_i$ is the number of times operator $i$ has been used, $s_i$ is the performance score of the operator $i$ and $\rho$ is the evaporation coefficient that prevents the procedure from converging too fast on a sub-optimal solution. Accordingly, the probability of using each operator is computed as:

\begin{gather}
    P(i)=\frac{W_i}{\sum_{j\in\Omega^+} W_j}\quad\forall i\in\Omega^+ \label{equ30}\\
    P(i)=\frac{W_i}{\sum_{j\in\Omega} W_j}\quad\forall i\in\Omega \label{equ31}\\
    P(i)=\frac{W_i}{\sum_{j\in\Omega^-} W_j}\quad\forall i\in\Omega^- \label{equ32}
\end{gather}

The acceptance criterion of the ALNS procedure is aligned with the principle used in the simulated annealing (SA) heuristic, as introduced by \citet{kirkpatrick1983optimization}. According to this framework, the probability of replacing the current solution with a new solution that performs worse is calculated using the formula $P=e^{-\frac{\Delta f(S)}{T}}$, where $\Delta f(S)$ represents the difference in net profits between the new worse solution and the current solution and $T$ is the current temperature in the SA process. The temperature decreases from an initial value $T_{init}$ over time based on a cooling rate $c$ which makes the ALNS procedure less likely to accept worse solutions as it progresses and encourages convergence towards an optimal solution. Next, the reheating step where the temperature $T$ is reset to its initial value $T_{init}$ when it falls to a minimum threshold $T_{min}$, ensures that the ALNS procedure retains some exploratory capability and avoids being trapped in local optima \citep{ali2021models}.

\subsection{A hybrid Dual-TS-ANS: the local search method}
In this section, we present the other crucial integral part of our MA - a hybrid Dual Tabu Search with Adaptive Neighborhood Search (Dual-TS-ANS) procedure. The Dual-TS-ANS functions as the local search mechanism, playing a role in the exploitation of the search space by systematically improving upon the current solution while avoiding cycles and preventing the MA from getting trapped in local optima. Our Dual-TS-ANS is inspired by the tabu search procedure in \citet{luMemeticAlgorithmOrienteering2018}, but has been adapted and refined to address the unique requirements and challenges of the HHC\&HCRRP-RNC. The Dual-TS-ANS incorporates three components that significantly contribute to its effectiveness and flexibility. Firstly, the ANS part involves iterative processes where insertion, internal, and removal operators, identical to those in ALNS (Section \ref{ALNS}), are dynamically selected based on their weighted probabilities. Next, the TS part evaluates the new solution and applies the tabu criteria after each operation of ANS. This integration ensures that the search process benefits from both the diversified exploration of TS and the specialized operations of ANS. Notably, as its name suggests, we incorporate two distinct tabu lists in the Dual-TS along with their respective criteria. These are the tabu-node list and the tabu-solution list. Lastly but notably, we develop a violation mechanism that contributes to handling infeasible solutions. By adopting this approach, we can aim to explore the possibility of transforming these infeasible solutions into feasible ones, thereby expanding the search horizon.  

The overall structure of the TS-ANS is outlined in Algorithm \ref{alg3}, with its key elements described in the following sections.

\begin{algorithm}[h]
	\caption{Dual-tabu-search with adaptive neighborhood search}
	\label{alg3}
	\begin{algorithmic}[1]
		\REQUIRE \emph{$G=(V, E)$}: a weighted undirected graph; \emph{$S_0$}: a solution from the initial population generated by ALNS; \emph{$\Omega^+$}: the set of insertion operators; \emph{$\Omega$}: the set of internal move operators; \emph{$\Omega^-$}: the set of removal operators; \emph{$\Omega^v$}: the set of violation operators; \emph{$L_{TS}$}: search depth; \emph{$u_p$}: frequency of updating the penalty parameter \emph{$\phi$}; \emph{$s_1$}: the score of finding a better feasible solution; \emph{$s_2$}: the score of finding a better conditional solution under violation criterion; \emph{$s_3$}: the score of finding an unfeasible solution; \emph{$u_w$}: frequency of updating weights of operators in \emph{$\Omega^+,\Omega,\Omega^-$}
  		\ENSURE \emph{$S^*$}: best feasible solution found so far
		\STATE \emph{Iter} $\leftarrow0$ \quad/*iteration counter*/
        \STATE $S^*\leftarrow S_0$ \quad/*the best feasible solution*/
        \STATE $S'\leftarrow S_0$ \quad/*solution with the best objective value*/
        \STATE Initialize tabu-node list, tabu-solution list, and their respective tabu tenure
        \WHILE{\emph{Iter} < $L_{TS}$}
        \STATE Select an operator \emph{$\omega^+,\omega$} and \emph{$\omega^-$} sequentially from \emph{$\Omega^+,\Omega$} and \emph{$\Omega^-$}
        \STATE $S\leftarrow$ the solution with a better objective function 
        \STATE $S \leftarrow \omega^+,\omega,\omega^-(S)$
        \IF{\emph{$S$} is a feasible solution and \emph{$f(S)>f(S^*)$}}
        \STATE \emph{$S^*\leftarrow S$}
        \STATE \emph{$Score(\omega^+,\omega,\omega^-)\;+=\;s_1$}
        \STATE update tabu-node list
        \ELSIF{\emph{$S$} is a non-tabu solution and does not involve tabu nodes}
        \IF{\emph{$S$} is conditionally feasible under violation criterion and \emph{$\Phi(S)>\Phi(S')$}}
        \STATE \emph{$S'\leftarrow S$}
        \STATE \emph{$Score(\omega^+,\omega,\omega^-)\;+=\;s_2$}
        \STATE \emph{$Score(\omega^v)\;+=\;s_2$}
        \STATE update tabu-node list
        \ELSE
        \STATE \emph{$Score(\omega^+,\omega,\omega^-)\;+=\;s_3$}
        \STATE \emph{$Score(\omega^v)\;+=\;s_3$}
        \STATE update tabu-solution list
        \ENDIF
        \ENDIF   
        \IF{\emph{Iter} mod \emph{$u_w$} = 0}
        \STATE Update Weight($\Omega^+,\Omega,\Omega^-,\Omega^v$)
        \STATE Score($\Omega^+,\Omega,\Omega^-,\Omega^v)\leftarrow0$
        \STATE Times($\Omega^+,\Omega,\Omega^-,\Omega^v)\leftarrow0$
        \ENDIF
        \IF{\emph{Iter} mod \emph{$u_p$} = 0}
        \STATE Update Weight($\Omega^v$)
        \STATE Score($\Omega^v)\leftarrow0$
        \STATE Times($\Omega^v)\leftarrow0$
        \ENDIF
        \ENDWHILE
	\end{algorithmic}
\end{algorithm}

\subsubsection{Violation mechanism and evaluation functions}
In every iteration, the Dual-TS-ANS dynamically selects between the two violation operators $\omega^v_1$ and $\omega^v_2$, which are distinct in their functions: when $\omega^v_1$ is chosen, it allows for the violation of the total travel time constraint, while the other operator permits the breach of the disruption time constraint. Thus, the solution is deemed conditionally feasible if it meets the disruption time constraint under the application of $\omega^v_1$; similarly, a solution is also considered conditionally feasible if it adheres to the total travel time constraint when $\omega^v_2$ is in effect. The degree to which the two constraints can be violated is regulated by a self-adjusting penalty parameter $\phi$ which modulates the level of infeasibility permitted in the search process, allowing our algorithm to explore solutions that initially violate certain constraints but potentially yield superior solutions in the following iterations. Inspired by the evaluation function design in \citet{luMemeticAlgorithmOrienteering2018}, the two evaluation functions under violation criteria are defined as follows:
\begin{gather}
    \Phi_1(S)=f(S)-\phi\times EX_1(S)\\
    \Phi_2(S)=f(S)-\phi\times EX_2(S)
\end{gather}
where $EX_1(S)=\max\{t(S)-T_{max},0\}$ represents the violation level of the total travel time budget, indicating how much time of solution $S$ surpasses the upper limit $T_{max}$. $EX_2(S)$ sums up the payment of existing customers who are rescheduled to experience disruptions of arrival times that exceed the maximum threshold. It reflects the extent of violation regarding the disruption time constraint. In addition, the penalty coefficient $\phi$ is initially set at 1 and is dynamically modified based on aims to maintain a balance in searching both feasible and infeasible domains. If the algorithm over-focuses on feasible solutions, we prefer $\phi$ with a lower value to encourage a deeper exploration into infeasible neighborhoods. Conversely, encountering too many infeasible solutions triggers an increase in $\phi$, steering the search back toward feasibility. Similar to the ALNS procedure, we also devise a selection mechanism for violation operators. The probabilities of using $\omega^v_1$ and $\omega^v_2$ are calculated following a method analogous to that outlined in formula (\ref{equ30}) - (\ref{equ32}).

\subsubsection{Dual tabu list and aspiration criterion}
Our proposed algorithm employs two types of tabu lists, namely tabu-node list and tabu-solution list. The tabu-node list is designed to prevent cycling actions in the following within the tabu tenure. It achieves this by prohibiting the revisitation of new customer nodes that were recently removed from the current solution, as well as preventing the removal of new customer nodes that were recently added, for the next $tl$ iterations (tabu-node tenure). On the other hand, the tabu-solution list is specifically created to exclude solutions that fail to meet conditionally feasible standards under violation criteria. Essentially, this means solutions that satisfy neither the total travel time constraint nor the disruption time constraint are deemed unsuitable for the subsequent iterations. We hold the view that such solutions are significantly less likely to contribute to the creation of superior feasible solutions in subsequent actions. 

The tabu tenure of tabu-node list, denoted as $tl(node)$, is dynamically set as $tl(node)=\alpha\times|N|+random(10)$, where $random(10)$ generates a random number between 1 and 10, $|N|$ represents the size of the new customer nodes, and $\alpha$ is the coefficient that controls tabu-node tenure.  For the tabu-solution list, its tenure is set the same as the population size, since it can help ensure that the algorithm avoids revisiting recent solutions prematurely, thereby encouraging exploration in diverse areas of the solution space. The aspiration criterion is utilized solely in conjunction with the tabu-node list. It is activated if a move results in an improved feasible solution, allowing for an exception to the tabu status.

\subsection{Population management}
The most fundamental method for updating the population typically involves a fitness-based replacement strategy, which replaces the least fit individuals in the population with new improved ones, to encourage the retention of high-quality solutions. Yet this method can cause premature convergence since the population tends to be composed of high-fitness individuals that are highly similar to each other. Therefore, our proposed MA implements a hybrid of elitist and diversity-preserving strategies, which strikes a balance between retaining high-quality solutions (elitists) and exploring new solution spaces (diversity preservation). Firstly, we divide the population into two tiers: elitist and average. The elitist tier contains the top 20\% of the population in terms of fitness, while the average tier encompasses the rest. Secondly, we deliver different replacement strategies for the two tiers of individuals. For the elitist tier, a simple fitness-based replacement strategy is employed where only better solutions replace the existing ones. In contrast, for the average tier, a fitness-distance ratio (FDR) strategy is adopted that both considers the fitness of one solution which evaluates the effectiveness of our proposed MA, and the dissimilarities among solutions which represent the diversity within the population. In other words, the FDR strategy prefers solutions that are not only high-performing but also different from the current population. For example, when choosing between two solutions with the same fitness, the one that is more different from the rest of the population will be selected. The FDR is computed as $FDR(S_i)=\frac{f(S_i)}{D(S_i)+\epsilon}$, where $D(S_i)$ is the average distance of solution $S_i$ to all other solutions and $\epsilon$ is a small constant to avoid division by zero. Mathematically, $D(S_i)$ is calculated by Equation (\ref{equ35}).
\begin{equation}
\label{equ35}
    D(S_i)=\frac{1}{n}\sum_{j=1,j\neq i}^n d(S_i,S_j)\quad\forall i\in P
\end{equation}
where $d(S_i,S_j)$ is the distance between solutions $S_i$ and $S_j$ within Population $P=\{S_i,\cdots,S_n\}$.

\section{Computational experiments} \label{sec_experiment}
In this section, we conduct detailed experiments to assess the performance of the proposed algorithm to address the HHC\&HCRRP-RNC. Note that the HHC\&HCRRP-RNC is a special case of the OPMV,  thus we use the benchmark instances provided by \citet{chao1993new,chao1996team}, which are frequently used in testing OP and the team OP (TOP), as the input for performance evaluation over all the stand-alone scenarios. Our adapted instances are grouped into three problem sizes based on the total count of existing and new customer nodes, i.e. small instances with no more than 25 nodes, medium instances with up to 50 nodes, and large instances with a maximum of 100 nodes. The evaluation involves three key components: first, we determine the parameters involved in the algorithm. Second, we assess the effectiveness of the proposed MA (MA2) against four approaches: MIP model, ALNS, TS, and MA1 \citep{luMemeticAlgorithmOrienteering2018} on small and medium-sized instances. Third, we compare the performance of ALNS, TS, MA1, and MA2 across all available instances. In our study, the MILP model is solved with Gurobi 11.0.0, and the four algorithms are coded in Python and performed on an Apple M1 pro chip (8-core CPU) with 16.0 GB  RAM running the macOS operating system. The runtime of all experiments for each instance is limited by 1800 CPU seconds and all algorithms are run 10 times.

\subsection{Parameter tuning}\label{sec_paras}
\subsubsection{Estimates for the total travel time and disruption time constraints}
In this section, we report a series of experiments conducted to determine the associated parameters. Firstly, in MA2 and the three comparative algorithms, we all apply two adaptive scaling parameters to discuss estimates for the total travel time and disruption time constraints in the HHC\&HCRRP-RNC: parameter $\mu$ controls the threshold for the maximum total travel time in the HHC\&HCRRP-RNC, and parameter $\lambda$ decides the maximum allowed disruption time.  

More precisely, to express the maximum allowed total travel time, we can construct Equation (\ref{equ36}):
\begin{equation}
\label{equ36}
    T_{max}=\mu\times T_{max}(TSP)
\end{equation}
where $T_{max}(TSP)$ represents a baseline travel time estimate which is derived from solving a TSP considering all existing and new customer nodes. For example, if a case involves 15 existing customer nodes and 10 new customer nodes, the baseline travel time limit is set the same as the optimal or near-optimal outcome (this depends on the feasibility of calculating the optimal solution within the given time limit) of the TSP that encompasses these 25 nodes.  The TSP outcome is solved using Gurobi solver, subject to a time limit of 1800 seconds. The parameter $\mu$ is applied along with the baseline travel time to establish an upper limit for the total travel time within the rescheduling context. It serves as a multiplier that scales the TSP-derived travel time to an appropriate level to provide flexibility to adjust to realistic operational constraints. For instance, a higher value of $\mu$ would result in a more lenient maximum travel time, allowing for longer routes in the rescheduled plan. However, an excessively high $\mu$ value can lead to impractically long routes that are not manageable within a single work period. Also, a high value could reduce the incentive for the algorithm to find more efficient and shorter routes since the extended time limit makes it easier to meet the travel time constraints without needing to optimize. In contrast, a smaller value of $\mu$ would impose a tighter limit to ensure more compact routing solutions, however, if $\mu$ is too low, it might render many potential routes infeasible, thereby making it more challenging to optimally utilize the resources. 

Moreover, we propose a simple model based on the average shortest path length to estimate the maximum allowed disruption time. Compared to other graph attributes summarized by \citet{liu2021time}, the average shortest path length considers the shortest paths between all pairs of nodes, offering a robust measure that is less affected by outliers or extreme values and being able to indicate how efficiently resources can be transferred in the graph. Combined with the parameter $\lambda$, the maximum allowed disruption time is computed as:
\begin{equation}
    T=\lambda\times\frac{1}{n\times(n-1)}\sum_{i\in V}\sum_{j\in V,j\neq i}d(i,j)
\end{equation}
where the fractions that follow $\lambda$ in the formula represent the average shortest path length in the graph, $d(i,j)$ represents the shortest path distance between vertices i and j, and $n$ is the number of vertices in $G$.

We test 10 and three different values for the parameter $\mu$ and $\lambda$ respectively, resulting in a total of 30 pairs of parameter combinations. Table \ref{tab2} reports the instances classified by the problem size, and Table \ref{tab3} lists the different settings of the parameter $\mu$ and $\lambda$.

\begin{table}[htbp]
\centering
	\setlength{\abovecaptionskip}{0pt}
	\setlength{\belowcaptionskip}{5pt}
	\caption{\normalsize{Benchmark instances}}
 	\begin{tabular}{cccc}
        \hline
\multirow{2}{*}{Problem size} & \multirow{2}{*}{\begin{tabular}[c]{@{}c@{}}Number \\of instances \end{tabular}}& \multirow{2}{*}{\begin{tabular}[c]{@{}c@{}} Existing \\customers \end{tabular}} & \multirow{2}{*}{\begin{tabular}[c]{@{}c@{}} New \\customers \end{tabular}}\\ 
&&&\\
\hline
\multirow{3}{*}{Small size} & 30 & 5 & 5\\
& 30 & 10 & 10\\
& 30 & 15 & 10\\ 
\hline
\multirow{3}{*}{Medium size} & 30  & 15 & 15\\
& 30 & 20 & 20\\
& 30 & 30 & 20\\ 
\hline
\multirow{3}{*}{Large size} & 30 & 30 & 30\\
& 30 & 40& 40\\
& 30 & 50 & 50\\
\hline
\end{tabular}
\label{tab2}
\end{table}

\begin{table}[htbp]
\centering
	\setlength{\abovecaptionskip}{0pt}
	\setlength{\belowcaptionskip}{5pt}
	\caption{\normalsize{Parameters summary in the HHC\&HCRRP-RNC}}
    \begin{tabular}{p{1.5cm}<{\centering\arraybackslash}p{8cm}<{\raggedright}p{3cm}<{\raggedright}}
\specialrule{0.05em}{3pt}{3pt}
Parameter & Description & Values\\\specialrule{0.05em}{3pt}{3pt}
\multirow{2}{*}{$\mu$} & Scaling factor that calculates the maximum allowed total travel time in the rescheduling problem & \multirow{2}{*}{0.8, 0.9, 1}\\\specialrule{0em}{1pt}{1pt}
\multirow{2}{*}{$\lambda$} & scaling factor that calculates the maximum allowed disruption time in the rescheduling problem & 0.1, 0.2, 0.3, 0.4, 0.5, 0.6, 0.7, 0.8, 0.9, 1\\\specialrule{0.05em}{2pt}{0pt}
\end{tabular}
\label{tab3}
\end{table}

\subsubsection{MA parameters tuning}
In this section, we determine the appropriate parameters associated with MA2 mentioned in Section \ref{sec_algorithm} for solving the HHC\&HCRRP-RNC. For enhanced clarity in this paper, we summarize these parameters and their corresponding tested values in Table \ref{tab4}.
\begin{table}[htbp]
\centering
	\renewcommand\arraystretch{1.5}
	\setlength{\abovecaptionskip}{0pt}
	\setlength{\belowcaptionskip}{5pt}
	\caption{\normalsize{Results from Parameter tuning of the algorithms}}
\begin{tabular}{p{1.5cm}<{\centering\arraybackslash}p{2cm}<{\centering\arraybackslash}p{5cm}<{\raggedright}p{3cm}<{\raggedright}p{1.5cm}<{\raggedright}}
\specialrule{0.05em}{3pt}{3pt}
Parameter & Algorithm & Description & Considered values & Final value\\
\specialrule{0.05em}{3pt}{3pt}
$N$ & Algorithm 1 & Population size & \{5, 10, 20\} & 5\\
$G$ & Algorithm 1 & The maximum generations & \{5, 10, 20\} & 5\\
$L_{ALNS}$ & Algorithm 2 & Search depth of ALNS & \{500, 1000, 1500\} & 1000\\
\multirow{2}{*}{$u_w$} & Algorithm 2, Algorithm 3 & Frequency of updating weights of operators in $\Omega^+$, $\Omega$, $\Omega^-$ and $\Omega^v$ & \multirow{2}{*}{\{50, 100, 200\}} & \multirow{2}{*}{50}\\
$\rho$ & Algorithm 2 & The evaporation factor & \{0.8, 0.9, 0.95\}  & 0.95\\
$u$ & Algorithm 2 & The perturbation factor & \{0.05, 0.1, 0.2\}  & 0.1\\
$L_{TS}$ & Algorithm 3  & Search depth of TS  & \{100, 150, 200\} & 200\\
\multirow{2}{*}{$u_p$} & \multirow{2}{*}{Algorithm 3} & Frequency of updating the penalty parameter $\phi$ & \multirow{2}{*}{\{5, 10, 15\}} & \multirow{2}{*}{10}\\
\specialrule{0.05em}{3pt}{3pt}
\end{tabular}
\label{tab4}
\end{table}

\subsection{Experiment 1: Comparisons with exact approaches on small and medium instances}
To validate the effectiveness of our proposed algorithm, this section compares the results of MA2 with those solved with Gurobi solver 11.0.0 to optimality, or near-optimality within the predefined time limit, according to the proposed mathematical model in Section \ref{sec_algorithm}.  However, due to the inability of Gurobi solver to handle very large size NP-hard problem instances competently, we pick the 180 instances of small and medium problem size(see Table \ref{tab1}) for the performance test in this section. In addition, for each instance, we implement algorithm MA2 10 times, record the best, average, and worst objective values that we obtained, and also record the best CPU time (see columns 8-14 in Table \ref{tab5}). For each instance, the gap is computed as follows:
\begin{equation}
    Gap=\frac{Gurobi\;objecive-MA2\;objective}{Gurobi\;objecive}\times100
\end{equation}

\begin{landscape}
\label{tab5}
\begin{longtable}{ccccllllllllll}
\caption{Results from Parameter tuning of the algorithms} \label{tab:my_label} \\
\toprule
$V^E$ & $V^O$ & $\lambda$ & $\mu$ & \multicolumn{3}{c}{Gurobi} & \multicolumn{7}{c}{MA2} \\
\cmidrule(lr){5-7} \cmidrule(lr){8-14}
 & & & & Obj. & Gap (\%) & Time (s) & Obj. (best) & Gap (\%) & Obj. (avg.) & Gap (\%) & Obj. (worst) & Gap (\%) & Time (s) \\
\midrule
\endfirsthead

\multicolumn{14}{c}%
{{\bfseries Table \thetable\ Continued from previous page}} \\
\toprule
$V^E$ & $V^O$ & $\lambda$ & $\mu$ & \multicolumn{3}{c}{Gurobi} & \multicolumn{7}{c}{MA2} \\
\cmidrule(lr){5-7} \cmidrule(lr){8-14}
 & & & & Obj. & Gap (\%) & Time (s) & Obj. (best) & Gap (\%) & Obj. (avg.) & Gap (\%) & Obj. (worst) & Gap (\%) & Time (s) \\
\midrule
\endhead

\midrule
\multicolumn{14}{r}{{Continued on next page}} \\
\endfoot

\bottomrule
\endlastfoot

5  & 5  & 1   & 1   & 126  & 0      & \textless{}1 & 114  & 9.52  & 103.2  & 15.24 & 102  & 19.05 & 2.19   \\
   &    &     & 0.9 & 102  & 0      & 2.25         & 102  & 0     & 91.2   & 10.59 & 90   & 11.76 & 2.88   \\
   &    &     & 0.8 & 72   & 0      & 2.22         & 72   & 0     & 72     & 0     & 72   & 0     & 1.99   \\
   &    & 0.9 & 1   & 126  & 0      & \textless{}1 & 108  & 9.52  & 103.2  & 15.24 & 102  & 19.05 & 3.25   \\
   &    &     & 0.9 & 102  & 0      & 1.94         & 96   & 5.88  & 91.8   & 10    & 90   & 11.76 & 3      \\
   &    &     & 0.8 & 72   & 0      & 2.03         & 72   & 0     & 72     & 0     & 72   & 0     & 1.86   \\
   &    & 0.8 & 1   & 126  & 0      & \textless{}1 & 114  & 9.52  & 104.4  & 17.14 & 102  & 19.05 & 2.16   \\
   &    &     & 0.9 & 102  & 0      & 3.31         & 102  & 0     & 102    & 0     & 102  & 0     & 1.93   \\
   &    &     & 0.8 & 72   & 0      & 2.04         & 72   & 0     & 72     & 0     & 72   & 0     & 1.88   \\
   &    & 0.7 & 1   & 108  & 0      & \textless{}1 & 108  & 0     & 102.6  & 5     & 102  & 5.56  & 3.65   \\
   &    &     & 0.9 & 90   & 0      & 1.43         & 90   & 0     & 85.2   & 5.33  & 84   & 6.67  & 2.01   \\
   &    &     & 0.8 & 72   & 0      & 1.04         & 72   & 0     & 72     & 0     & 72   & 0     & 2.42   \\
   &    & 0.6 & 1   & 108  & 0      & \textless{}1 & 108  & 0     & 105    & 2.78  & 102  & 5.56  & 2.76   \\
   &    &     & 0.9 & 90   & 0      & 1.21         & 90   & 0     & 85.2   & 5.33  & 84   & 6.67  & 2.55   \\
   &    &     & 0.8 & 72   & 0      & \textless{}1 & 72   & 0     & 72     & 0     & 72   & 0     & 1.99   \\
   &    & 0.5 & 1   & 108  & 0      & \textless{}1 & 108  & 0     & 106.8  & 1.11  & 102  & 5.56  & 2.06   \\
   &    &     & 0.9 & 90   & 0      & \textless{}1 & 90   & 0     & 88.8   & 1.33  & 84   & 6.67  & 2.68   \\
   &    &     & 0.8 & 72   & 0      & \textless{}1 & 72   & 0     & 72     & 0     & 72   & 0     & 1.81   \\
   &    & 0.4 & 1   & 108  & 0      & \textless{}1 & 108  & 0     & 108    & 0     & 108  & 0     & 1.88   \\
   &    &     & 0.9 & 90   & 0      & \textless{}1 & 90   & 0     & 90     & 0     & 90   & 0     & 2.03   \\
   &    &     & 0.8 & 72   & 0      & \textless{}1 & 72   & 0     & 72     & 0     & 72   & 0     & 1.88   \\
   &    & 0.3 & 1   & 108  & 0      & \textless{}1 & 108  & 0     & 105.6  & 2.22  & 96   & 11.11 & 2.18   \\
   &    &     & 0.9 & 84   & 0      & \textless{}1 & 84   & 0     & 84     & 0     & 84   & 0     & 1.88   \\
   &    &     & 0.8 & 72   & 0      & \textless{}1 & 72   & 0     & 72     & 0     & 72   & 0     & 1.88   \\
   &    & 0.2 & 1   & 108  & 0      & \textless{}1 & 108  & 0     & 108    & 0     & 108  & 0     & 2.08   \\
   &    &     & 0.9 & 84   & 0      & \textless{}1 & 84   & 0     & 84     & 0     & 84   & 0     & 2.09   \\
   &    &     & 0.8 & 72   & 0      & \textless{}1 & 72   & 0     & 72     & 0     & 72   & 0     & 2.48   \\
   &    & 0.1 & 1   & 108  & 0      & \textless{}1 & 108  & 0     & 108    & 0     & 108  & 0     & 2.03   \\
   &    &     & 0.9 & 84   & 0      & \textless{}1 & 84   & 0     & 84     & 0     & 84   & 0     & 2.03   \\
   &    &     & 0.8 & 72   & 0      & \textless{}1 & 72   & 0     & 72     & 0     & 72   & 0     & 2.3    \\
\hline
10 & 10 & 1   & 1   & 324  & 0      & 42.03        & 312  & 3.7   & 312    & 3.7   & 312  & 3.7   & 10.69  \\
   &    &     & 0.9 & 312  & 0      & 3.47         & 288  & 7.69  & 277.2  & 11.15 & 264  & 15.38 & 6.82   \\
   &    &     & 0.8 & 270  & 0      & 1.38         & 264  & 2.22  & 257.4  & 4.67  & 246  & 8.89  & 11.13  \\
   &    & 0.9 & 1   & 324  & 0      & 27.89        & 312  & 3.7   & 307.2  & 5.19  & 288  & 11.11 & 7.72   \\
   &    &     & 0.9 & 288  & 0      & 25.73        & 288  & 0     & 271.8  & 5.63  & 264  & 8.33  & 7.64   \\
   &    &     & 0.8 & 264  & 0      & 3.69         & 264  & 0     & 260.4  & 1.36  & 252  & 4.55  & 7.14   \\
   &    & 0.8 & 1   & 312  & 0      & 21.71        & 324  & 0     & 312    & 3.7   & 306  & 5.56  & 8.86   \\
   &    &     & 0.9 & 288  & 0      & 14.19        & 288  & 0     & 272.4  & 5.42  & 264  & 8.33  & 6.67   \\
   &    &     & 0.8 & 264  & 0      & 1.1          & 264  & 0     & 258    & 2.27  & 252  & 4.55  & 6.58   \\
   &    & 0.7 & 1   & 324  & 0      & 27.71        & 312  & 3.7   & 307.2  & 5.19  & 288  & 11.11 & 9.12   \\
   &    &     & 0.9 & 288  & 0      & 12.71        & 270  & 6.25  & 264.6  & 8.12  & 264  & 8.33  & 11.2   \\
   &    &     & 0.8 & 264  & 0      & \textless{}1 & 264  & 0     & 262.8  & 0.45  & 252  & 4.55  & 13.16  \\
   &    & 0.6 & 1   & 324  & 0      & 12.3         & 318  & 1.85  & 312.6  & 3.52  & 312  & 3.7   & 13.95  \\
   &    &     & 0.9 & 288  & 0      & 6.29         & 282  & 2.08  & 277.2  & 3.75  & 264  & 8.33  & 13.15  \\
   &    &     & 0.8 & 264  & 0      & \textless{}1 & 264  & 0     & 260.4  & 1.36  & 252  & 4.55  & 13.27  \\
   &    & 0.5 & 1   & 324  & 0      & 13.1         & 312  & 3.7   & 312    & 3.7   & 312  & 3.7   & 13.56  \\
   &    &     & 0.9 & 276  & 0      & 8.87         & 276  & 0     & 268.8  & 2.61  & 264  & 4.35  & 16.03  \\
   &    &     & 0.8 & 264  & 0      & \textless{}1 & 264  & 0     & 259.8  & 1.59  & 246  & 6.82  & 17.02  \\
   &    & 0.4 & 1   & 324  & 0      & 17.7         & 312  & 3.7   & 312    & 3.7   & 312  & 3.7   & 17.28  \\
   &    &     & 0.9 & 276  & 0      & 7.82         & 270  & 2.17  & 267.6  & 3.04  & 264  & 4.35  & 14.58  \\
   &    &     & 0.8 & 252  & 0      & 1.23         & 246  & 2.38  & 246    & 2.38  & 246  & 2.38  & 12.99  \\
   &    & 0.3 & 1   & 324  & 0      & 4.71         & 312  & 3.7   & 304.2  & 6.11  & 282  & 12.96 & 10.74  \\
   &    &     & 0.9 & 276  & 0      & 8.23         & 270  & 2.17  & 265.2  & 3.91  & 264  & 4.35  & 13.79  \\
   &    &     & 0.8 & 252  & 0      & \textless{}1 & 252  & 0     & 244.2  & 3.1   & 234  & 7.14  & 14.07  \\
   &    & 0.2 & 1   & 324  & 0      & 3.81         & 324  & 0     & 306.6  & 5.37  & 270  & 16.67 & 12.46  \\
   &    &     & 0.9 & 276  & 0      & 6.93         & 270  & 2.17  & 265.8  & 3.7   & 264  & 4.35  & 10.72  \\
   &    &     & 0.8 & 246  & 0      & \textless{}1 & 246  & 0     & 242.4  & 1.46  & 228  & 7.32  & 7.7    \\
   &    & 0.1 & 1   & 324  & 0      & 4.19         & 324  & 0     & 316.8  & 2.22  & 312  & 3.7   & 9.86   \\
   &    &     & 0.9 & 276  & 0      & 4.34         & 276  & 0     & 265.8  & 3.7   & 264  & 4.35  & 18.25  \\
   &    &     & 0.8 & 246  & 0      & \textless{}1 & 246  & 0     & 243    & 1.22  & 228  & 7.32  & 13.69  \\
\hline
15 & 10 & 1   & 1   & 450  & 0      & 4.7          & 450  & 0     & 442.2  & 1.73  & 420  & 6.67  & 14.23  \\
   &    &     & 0.9 & 408  & 0      & 2.17         & 408  & 0     & 394.2  & 3.38  & 384  & 5.88  & 15.97  \\
   &    &     & 0.8 & 324  & 0      & 4.09         & 324  & 0     & 324    & 0     & 324  & 0     & 15.92  \\
   &    & 0.9 & 1   & 450  & 0      & 5.73         & 450  & 0     & 445.2  & 1.07  & 426  & 5.33  & 17.77  \\
   &    &     & 0.9 & 408  & 0      & 3.19         & 408  & 0     & 389.4  & 4.56  & 372  & 8.82  & 17.77  \\
   &    &     & 0.8 & 324  & 0      & 2.48         & 324  & 0     & 321.6  & 0.74  & 318  & 1.85  & 11.64  \\
   &    & 0.8 & 1   & 450  & 0      & 3.17         & 450  & 0     & 429.6  & 4.53  & 420  & 6.67  & 14.3   \\
   &    &     & 0.9 & 408  & 0      & 1.42         & 402  & 1.47  & 398.4  & 2.35  & 384  & 5.88  & 14.06  \\
   &    &     & 0.8 & 324  & 0      & 3.52         & 324  & 0     & 319.2  & 1.48  & 318  & 1.85  & 15.31  \\
   &    & 0.7 & 1   & 450  & 0      & 11.45        & 444  & 1.33  & 429.6  & 4.53  & 420  & 6.67  & 12.16  \\
   &    &     & 0.9 & 408  & 0      & 2.67         & 408  & 0     & 402.6  & 1.32  & 402  & 1.47  & 18.12  \\
   &    &     & 0.8 & 324  & 0      & 3.43         & 324  & 0     & 319.8  & 1.3   & 318  & 1.85  & 13.72  \\
   &    & 0.6 & 1   & 444  & 0      & 12.11        & 444  & 0     & 429.6  & 3.24  & 420  & 5.41  & 17.42  \\
   &    &     & 0.9 & 402  & 0      & 1.33         & 402  & 0     & 398.4  & 0.9   & 384  & 4.48  & 16.17  \\
   &    &     & 0.8 & 324  & 0      & 2.3          & 324  & 0     & 321.6  & 0.74  & 312  & 3.7   & 18.73  \\
   &    & 0.5 & 1   & 444  & 0      & 8.4          & 444  & 0     & 420    & 5.41  & 408  & 8.11  & 19.22  \\
   &    &     & 0.9 & 402  & 0      & 1.14         & 402  & 0     & 394.8  & 1.79  & 360  & 10.45 & 13.43  \\
   &    &     & 0.8 & 324  & 0      & 1.12         & 324  & 0     & 321.6  & 0.74  & 318  & 1.85  & 11.59  \\
   &    & 0.4 & 1   & 444  & 0      & 10.39        & 444  & 0     & 424.8  & 4.32  & 408  & 8.11  & 17.65  \\
   &    &     & 0.9 & 402  & 0      & 2.21         & 402  & 0     & 391.8  & 2.54  & 378  & 5.97  & 20.89  \\
   &    &     & 0.8 & 324  & 0      & 1.08         & 324  & 0     & 323.4  & 0.19  & 318  & 1.85  & 13.72  \\
   &    & 0.3 & 1   & 444  & 0      & 1.21         & 444  & 0     & 429.6  & 3.24  & 408  & 8.11  & 18.94  \\
   &    &     & 0.9 & 402  & 0      & \textless{}1 & 402  & 0     & 385.2  & 4.18  & 366  & 8.96  & 19.61  \\
   &    &     & 0.8 & 324  & 0      & \textless{}1 & 324  & 0     & 320.4  & 1.11  & 318  & 1.85  & 15.52  \\
   &    & 0.2 & 1   & 420  & 0      & 4.56         & 420  & 0     & 408    & 2.86  & 396  & 5.71  & 18.19  \\
   &    &     & 0.9 & 402  & 0      & 1.03         & 396  & 1.49  & 368.4  & 8.36  & 360  & 10.45 & 14.97  \\
   &    &     & 0.8 & 324  & 0      & \textless{}1 & 324  & 0     & 322.2  & 0.56  & 318  & 1.85  & 16.84  \\
   &    & 0.1 & 1   & 420  & 0      & \textless{}1 & 420  & 0     & 409.8  & 2.43  & 378  & 10    & 19.92  \\
   &    &     & 0.9 & 378  & 0      & \textless{}1 & 378  & 0     & 373.8  & 1.11  & 360  & 4.76  & 16.49  \\
   &    &     & 0.8 & 324  & 0      & \textless{}1 & 324  & 0     & 323.4  & 0.19  & 318  & 1.85  & 12.63  \\
\hline
15 & 15 & 1   & 1   & 612  & 0      & 454.28       & 588  & 3.92  & 562.2  & 8.14  & 528  & 13.73 & 18.57  \\
   &    &     & 0.9 & 570  & 0      & 29.37        & 570  & 0     & 545.4  & 4.32  & 510  & 10.53 & 10.92  \\
   &    &     & 0.8 & 486  & 0      & 15           & 456  & 6.17  & 429    & 11.73 & 408  & 16.05 & 15.75  \\
   &    & 0.9 & 1   & 612  & 0      & 308.4        & 564  & 7.84  & 552.6  & 9.71  & 528  & 13.73 & 17.85  \\
   &    &     & 0.9 & 552  & 0      & 42.33        & 552  & 0     & 524.4  & 5     & 498  & 9.78  & 11.51  \\
   &    &     & 0.8 & 486  & 0      & 19.19        & 486  & 0     & 461.4  & 5.06  & 444  & 8.64  & 14.42  \\
   &    & 0.8 & 1   & 612  & 0      & 149.21       & 558  & 8.82  & 546.6  & 10.69 & 528  & 13.73 & 10.41  \\
   &    &     & 0.9 & 552  & 0      & 53.83        & 540  & 2.17  & 522.6  & 5.33  & 492  & 10.87 & 14.39  \\
   &    &     & 0.8 & 486  & 0      & 8.38         & 486  & 0     & 461.4  & 5.06  & 420  & 13.58 & 15.76  \\
   &    & 0.7 & 1   & 588  & 1.0204 & 1800         & 540  & 8.16  & 518.4  & 11.84 & 504  & 14.29 & 15.02  \\
   &    &     & 0.9 & 552  & 0      & 25.11        & 534  & 3.26  & 514.2  & 6.85  & 504  & 8.7   & 16.4   \\
   &    &     & 0.8 & 486  & 0      & 6.94         & 480  & 1.23  & 467.4  & 3.83  & 432  & 11.11 & 19.61  \\
   &    & 0.6 & 1   & 588  & 0      & 1734.99      & 558  & 5.1   & 543.6  & 7.55  & 504  & 14.29 & 13.9   \\
   &    &     & 0.9 & 540  & 0      & 63.95        & 528  & 2.22  & 520.2  & 3.67  & 504  & 6.67  & 10.24  \\
   &    &     & 0.8 & 486  & 0      & 3.69         & 486  & 0     & 479.4  & 1.36  & 450  & 7.41  & 12.82  \\
   &    & 0.5 & 1   & 588  & 0      & 320.7        & 552  & 6.12  & 531    & 9.69  & 510  & 13.27 & 13.23  \\
   &    &     & 0.9 & 528  & 0      & 42.51        & 528  & 0     & 518.4  & 1.82  & 504  & 4.55  & 11.46  \\
   &    &     & 0.8 & 486  & 0      & 7.49         & 480  & 1.23  & 471    & 3.09  & 456  & 6.17  & 14.94  \\
   &    & 0.4 & 1   & 588  & 0      & 407.07       & 582  & 1.02  & 544.8  & 7.35  & 510  & 13.27 & 18.65  \\
   &    &     & 0.9 & 528  & 0      & 39.05        & 528  & 0     & 524.4  & 0.68  & 510  & 3.41  & 13.59  \\
   &    &     & 0.8 & 486  & 0      & 2.98         & 486  & 0     & 472.8  & 2.72  & 462  & 4.94  & 11.93  \\
   &    & 0.3 & 1   & 588  & 0      & 100.47       & 558  & 5.1   & 547.2  & 6.94  & 528  & 10.2  & 12.15  \\
   &    &     & 0.9 & 528  & 0      & 24.33        & 528  & 0     & 504    & 4.55  & 474  & 10.23 & 15.19  \\
   &    &     & 0.8 & 486  & 0      & 1.43         & 468  & 3.7   & 450    & 7.41  & 414  & 14.81 & 18.5   \\
   &    & 0.2 & 1   & 588  & 0      & 50.19        & 582  & 1.02  & 537    & 8.67  & 510  & 13.27 & 13.14  \\
   &    &     & 0.9 & 528  & 0      & 21.19        & 516  & 2.27  & 487.2  & 7.73  & 456  & 13.64 & 10.92  \\
   &    &     & 0.8 & 486  & 0      & 1.92         & 486  & 0     & 463.8  & 4.57  & 438  & 9.88  & 12.91  \\
   &    & 0.1 & 1   & 588  & 0      & 24.1         & 582  & 1.02  & 555    & 5.61  & 510  & 13.27 & 12.05  \\
   &    &     & 0.9 & 528  & 0      & 16.95        & 516  & 2.27  & 489    & 7.39  & 456  & 13.64 & 14.7   \\
   &    &     & 0.8 & 468  & 0      & 5.43         & 468  & 0     & 448.8  & 4.1   & 426  & 8.97  & 36.95  \\
\hline
20 & 20 & 1   & 1   & 894  & 2.6846 & 1800         & 790  & 11.63 & 758.4  & 15.17 & 720  & 19.46 & 31.88  \\
   &    &     & 0.9 & 786  & 3.8168 & 1800         & 732  & 6.87  & 665.4  & 15.34 & 600  & 23.66 & 37.21  \\
   &    &     & 0.8 & 708  & 0      & 62.14        & 630  & 11.02 & 591.6  & 33.83 & 534  & 40.27 & 33.07  \\
   &    & 0.9 & 1   & 894  & 2.6846 & 1800         & 804  & 10.07 & 751.2  & 15.97 & 732  & 18.12 & 23.39  \\
   &    &     & 0.9 & 780  & 4.6154 & 1800         & 696  & 10.77 & 664.2  & 14.85 & 624  & 20    & 32.01  \\
   &    &     & 0.8 & 690  & 0      & 78.41        & 624  & 9.57  & 589.2  & 34.09 & 540  & 39.6  & 30.96  \\
   &    & 0.8 & 1   & 894  & 2.0134 & 1800         & 762  & 14.77 & 739.2  & 17.32 & 696  & 22.15 & 43.23  \\
   &    &     & 0.9 & 780  & 3.8462 & 1800         & 690  & 11.54 & 660.6  & 15.31 & 624  & 20    & 34.12  \\
   &    &     & 0.8 & 684  & 0      & 102.55       & 630  & 7.89  & 586.8  & 24.77 & 534  & 31.54 & 32.56  \\
   &    & 0.7 & 1   & 894  & 2.0134 & 1800         & 810  & 9.4   & 744.6  & 16.71 & 702  & 21.48 & 37.82  \\
   &    &     & 0.9 & 768  & 6.25   & 1800         & 732  & 4.69  & 665.4  & 13.36 & 618  & 19.53 & 29.93  \\
   &    &     & 0.8 & 684  & 0      & 60.32        & 624  & 8.77  & 600.6  & 12.19 & 582  & 14.91 & 25.56  \\
   &    & 0.6 & 1   & 894  & 1.3423 & 1800         & 798  & 10.74 & 755.4  & 1.64  & 684  & 10.94 & 40.53  \\
   &    &     & 0.9 & 768  & 4.6875 & 1800         & 726  & 5.47  & 673.8  & 12.27 & 636  & 17.19 & 37.62  \\
   &    &     & 0.8 & 684  & 0      & 77.19        & 642  & 6.14  & 597.6  & 12.63 & 582  & 14.91 & 21.48  \\
   &    & 0.5 & 1   & 894  & 0.6711 & 1800         & 786  & 12.08 & 749.4  & 16.17 & 708  & 20.81 & 38.78  \\
   &    &     & 0.9 & 762  & 4.7244 & 1800         & 720  & 5.51  & 673.2  & 1.58  & 636  & 7.02  & 41.42  \\
   &    &     & 0.8 & 678  & 0      & 53.6         & 624  & 7.96  & 608.4  & 10.27 & 594  & 12.39 & 26.13  \\
   &    & 0.4 & 1   & 894  & 0      & 678.99       & 798  & 10.74 & 762    & 14.77 & 738  & 17.45 & 33.27  \\
   &    &     & 0.9 & 762  & 3.937  & 1800         & 726  & 4.72  & 675    & 11.42 & 630  & 17.32 & 29.56  \\
   &    &     & 0.8 & 678  & 0      & 37.86        & 606  & 10.62 & 604.8  & 10.8  & 600  & 11.5  & 35.44  \\
   &    & 0.3 & 1   & 894  & 0      & 358.8        & 822  & 8.05  & 766.2  & 14.3  & 726  & 18.79 & 38.04  \\
   &    &     & 0.9 & 762  & 1.5748 & 1800         & 690  & 9.45  & 676.2  & 11.26 & 624  & 18.11 & 39.96  \\
   &    &     & 0.8 & 654  & 0      & 81.04        & 630  & 3.67  & 603.6  & 7.71  & 564  & 13.76 & 36.85  \\
   &    & 0.2 & 1   & 894  & 0      & 185.54       & 810  & 9.4   & 786.6  & 12.01 & 756  & 15.44 & 33.22  \\
   &    &     & 0.9 & 762  & 0.7874 & 1800         & 702  & 7.87  & 686.4  & 9.92  & 660  & 13.39 & 42.35  \\
   &    &     & 0.8 & 654  & 0      & 32.55        & 630  & 3.67  & 609    & 6.88  & 582  & 11.01 & 31.88  \\
   &    & 0.1 & 1   & 870  & 0      & 1195.56      & 822  & 5.52  & 774    & 11.03 & 732  & 15.86 & 38.87  \\
   &    &     & 0.9 & 756  & 0      & 1310.76      & 732  & 3.17  & 665.4  & 11.98 & 600  & 20.63 & 37.21  \\
   &    &     & 0.8 & 630  & 0      & 122.44       & 630  & 0     & 591.6  & 6.1   & 534  & 15.24 & 33.07  \\
\hline
30 & 20 & 1   & 1   & 1170 & 1.5385 & 1800         & 1104 & 5.64  & 1065.6 & 8.92  & 1020 & 12.82 & 97.88 \\
   &    &     & 0.9 & 1104 & 0      & 281.11       & 1002 & 9.24  & 979.2  & 11.3  & 960  & 13.04 & 72.83  \\
   &    &     & 0.8 & 996  & 0      & 89.78        & 900  & 9.64  & 873    & 12.35 & 846  & 15.06 & 78.21  \\
   &    & 0.9 & 1   & 1170 & 1.5385 & 1800         & 1170 & 0     & 1065.6 & 8.92  & 1020 & 12.82 & 65.08  \\
   &    &     & 0.9 & 1104 & 0      & 435.89       & 1044 & 5.43  & 1003.8 & 9.08  & 954  & 13.59 & 80.76  \\
   &    &     & 0.8 & 996  & 0      & 139.43       & 912  & 8.43  & 881.4  & 11.51 & 840  & 15.66 & 96.29  \\
   &    & 0.8 & 1   & 1170 &1.5385  &  1800        &1146  &2.05   &1084.2  &7.33   &1044  &10.77  &112.97 \\
   &    &     & 0.9 & 1086 & 1.105  & 1800         & 1050 & 3.31  & 997.8  & 8.12  & 954  & 12.15 & 114.49 \\
   &    &     & 0.8 & 972  & 0      & 315.34       & 936  & 3.7   & 883.2  & 9.14  & 852  & 12.35 & 89.2   \\
   &    & 0.7 & 1   & 1170 & 1.5385 & 1800         & 1128 & 3.59  & 1087.2 & 7.08  & 1026 & 12.31 & 77.82  \\
   &    &     & 0.9 & 1086 & 0      & 1598.76      & 1050 & 3.31  & 1011   & 6.91  & 960  & 11.6  & 66.35  \\
   &    &     & 0.8 & 972  & 0      & 177.93       & 960  & 1.23  & 894    & 8.02  & 864  & 11.11 & 62.76  \\
   &    & 0.6 & 1   & 1170 & 1.5385 & 1800         & 1122 & 4.1   & 1075.2 & 8.1   & 1026 & 12.31 & 117.54 \\
   &    &     & 0.9 & 1080 & 1.1111 & 1800         & 1020 & 5.56  & 977.4  & 9.5   & 936  & 13.33 & 75.05  \\
   &    &     & 0.8 & 960  & 0      & 410.53       & 918  & 4.38  & 888.6  & 7.44  & 864  & 10    & 74.71  \\
   &    & 0.5 & 1   & 1170 & 1.5385 & 1800         & 1092 & 6.67  & 1066.8 & 8.82  & 1026 & 12.31 & 97.41  \\
   &    &     & 0.9 & 1080 & 0.5556 & 1800         & 1008 & 6.67  & 981.6  & 9.11  & 918  & 15    & 98.88  \\
   &    &     & 0.8 & 960  & 0      & 188.01       & 918  & 4.38  & 869.4  & 9.44  & 840  & 12.5  & 66.6   \\
   &    & 0.4 & 1   & 1170 & 0      & 1708.13      & 1092 & 6.67  & 1054.8 & 9.85  & 1008 & 13.85 & 122.65 \\
   &    &     & 0.9 & 1080 & 0      & 1545.58      & 1008 & 6.67  & 966    & 10.56 & 930  & 13.89 & 99.31  \\
   &    &     & 0.8 & 960  & 0      & 107.12       & 882  & 8.13  & 861    & 10.31 & 840  & 12.5  & 84.2   \\
   &    & 0.3 & 1   & 1170 & 0      & 970.67       & 1140 & 2.56  & 1083   & 7.44  & 1050 & 10.26 & 121.94 \\
   &    &     & 0.9 & 1080 & 0      & 1350.62      & 1026 & 5     & 972.6  & 9.94  & 942  & 12.78 & 82.9   \\
   &    &     & 0.8 & 930  & 1.2903 & 1800         & 894  & 3.87  & 864.6  & 7.03  & 828  & 10.97 & 64.18  \\
   &    & 0.2 & 1   & 1170 & 0      & 268.89       & 1110 & 5.13  & 1066.2 & 8.87  & 1008 & 13.85 & 80.5   \\
   &    &     & 0.9 & 1080 & 0      & 246.16       & 1044 & 3.33  & 984.6  & 8.83  & 930  & 13.89 & 84.62  \\
   &    &     & 0.8 & 930  & 0      & 459.61       & 906  & 2.58  & 882.6  & 5.1   & 858  & 7.74  & 91.81  \\
   &    & 0.1 & 1   & 1170 & 0      & 117          & 1116 & 4.62  & 1071.6 & 8.41  & 1020 & 12.82 & 83.16  \\
   &    &     & 0.9 & 1062 & 0      & 1762.28      & 1002 & 5.65  & 972    & 8.47  & 942  & 11.3  & 95.72  \\
   &    &     & 0.8 & 906  & 0      & 851.55       & 900  & 0.66  & 880.2  & 2.85  & 858  & 5.3   & 69.16 
\end{longtable}
\end{landscape}

Notice in Table \ref{tab5} that Gurobi successfully solves 155 out of 180 instances (86.11\%) to optimality within the given allowed computational time. In contrast,  our proposed algorithm is able to find 81 optimal values, accounting for  45\% in the Obj. (best) column. Furthermore,  the Gap best reports that 72.22\% of the best solutions achieve a gap of 5\% or less. Also, the minor variations between Gap best, Gap avg., and Gap worst underscore the MA2 algorithm's robustness. The runtime for all instances remains below 100 CPU seconds, meanwhile, Gurobi's runtime exceeds 100 seconds in 59 instances, representing 32.78\% of the cases.

\subsection{Experiment 2: Comparisons with the existing approaches on all instances}
In this section, we present the results of the 90 large-size instances (see Table \ref{tab1}) obtained by the MA2 compared with those obtained by ALNS, TS, and MA1 \citep{luMemeticAlgorithmOrienteering2018}, in terms of the average solution quality.

\begin{figure}[htbp]
\centering
\begin{subfigure}{\textwidth}
  \centering
  \includegraphics[width=0.65\linewidth]{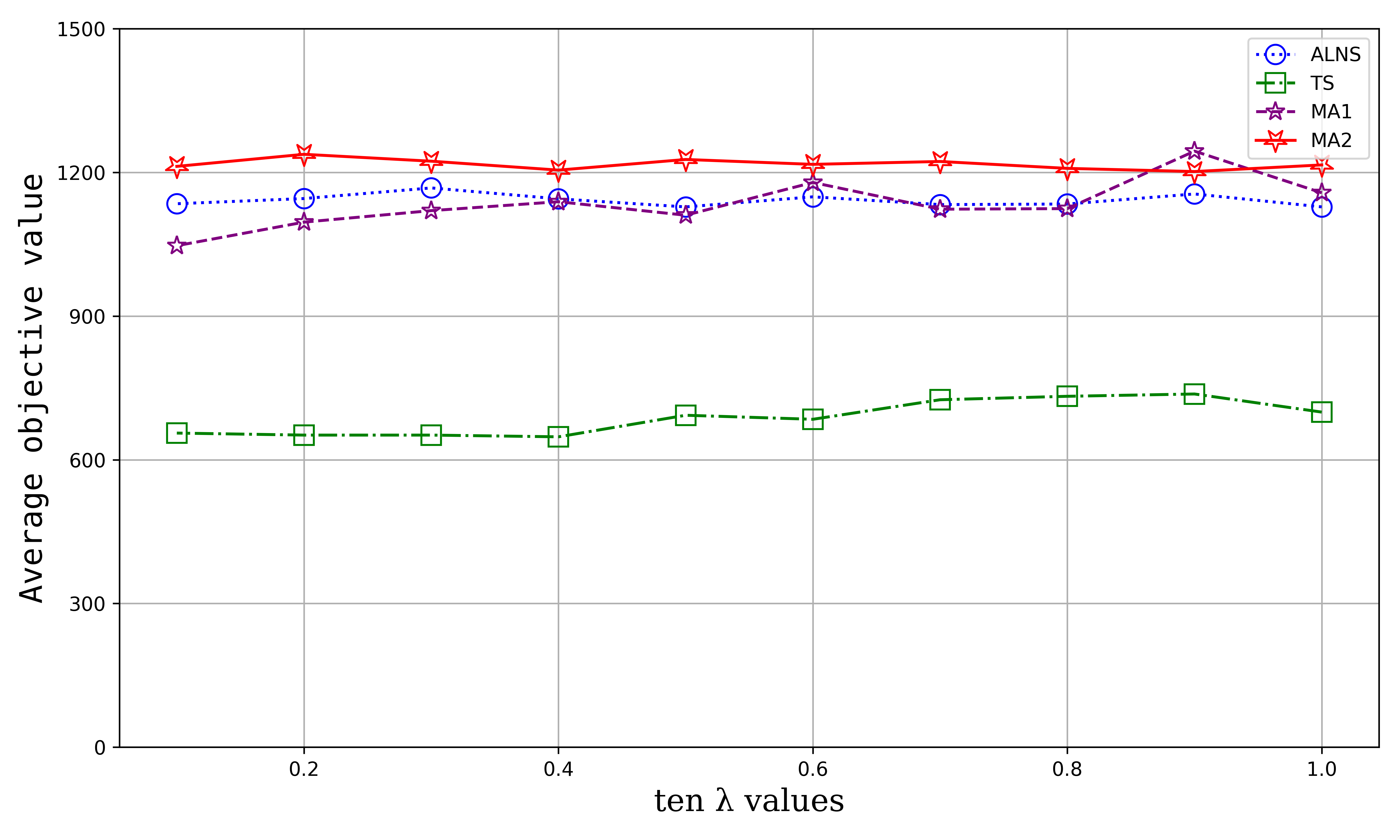}
  \caption{Existing customers=30, New customers=30, $\mu$=1}
\end{subfigure}

\begin{subfigure}{\textwidth}
  \centering
  \includegraphics[width=0.65\linewidth]{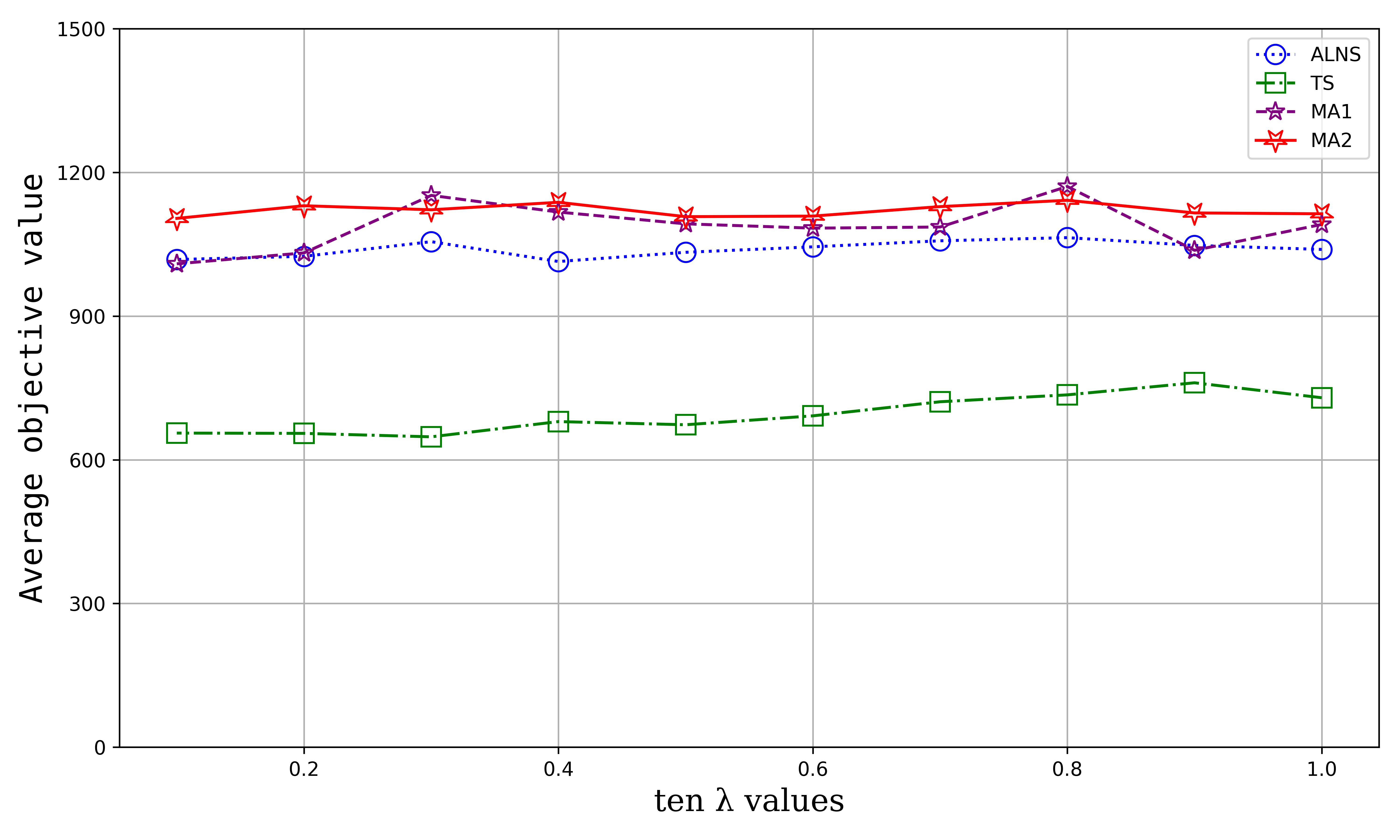}
  \caption{Existing customers=30, New customers=30, $\mu$=0.9}
\end{subfigure}

\begin{subfigure}{\textwidth}
  \centering
  \includegraphics[width=0.65\linewidth]{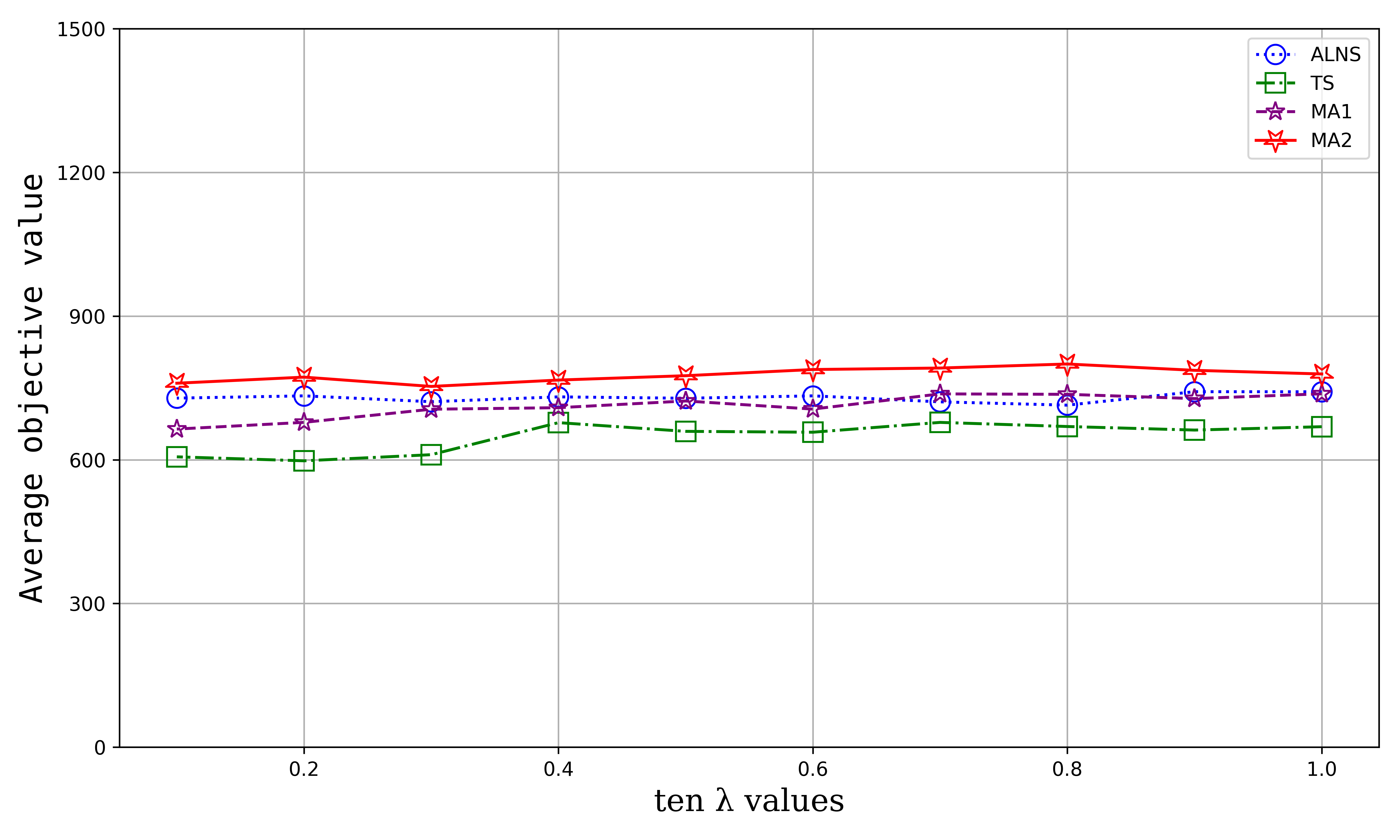}
  \caption{Existing customers=30, New customers=30, $\mu$=0.8}
\end{subfigure}
\caption{Average objective values obtained by MA2, ALNS, TS, and MA1 with 30 existing and 30 new customers}
\label{fig6}
\end{figure}

\begin{figure}[htbp]
\centering
\begin{subfigure}{\textwidth}
  \centering
  \includegraphics[width=0.65\linewidth]{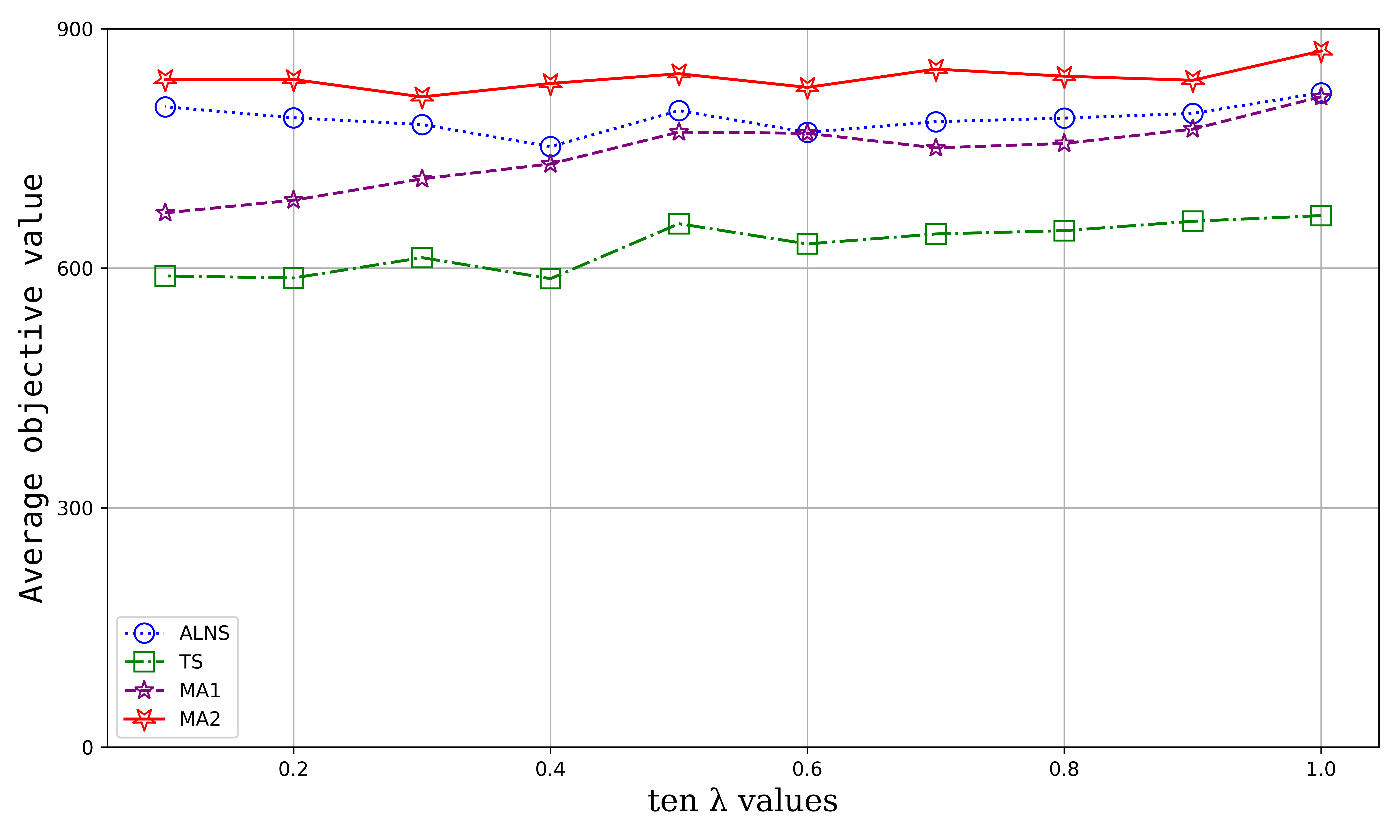}
  \caption{Existing customers=40, New customers=40, $\mu$=1}
\end{subfigure}

\begin{subfigure}{\textwidth}
  \centering
  \includegraphics[width=0.65\linewidth]{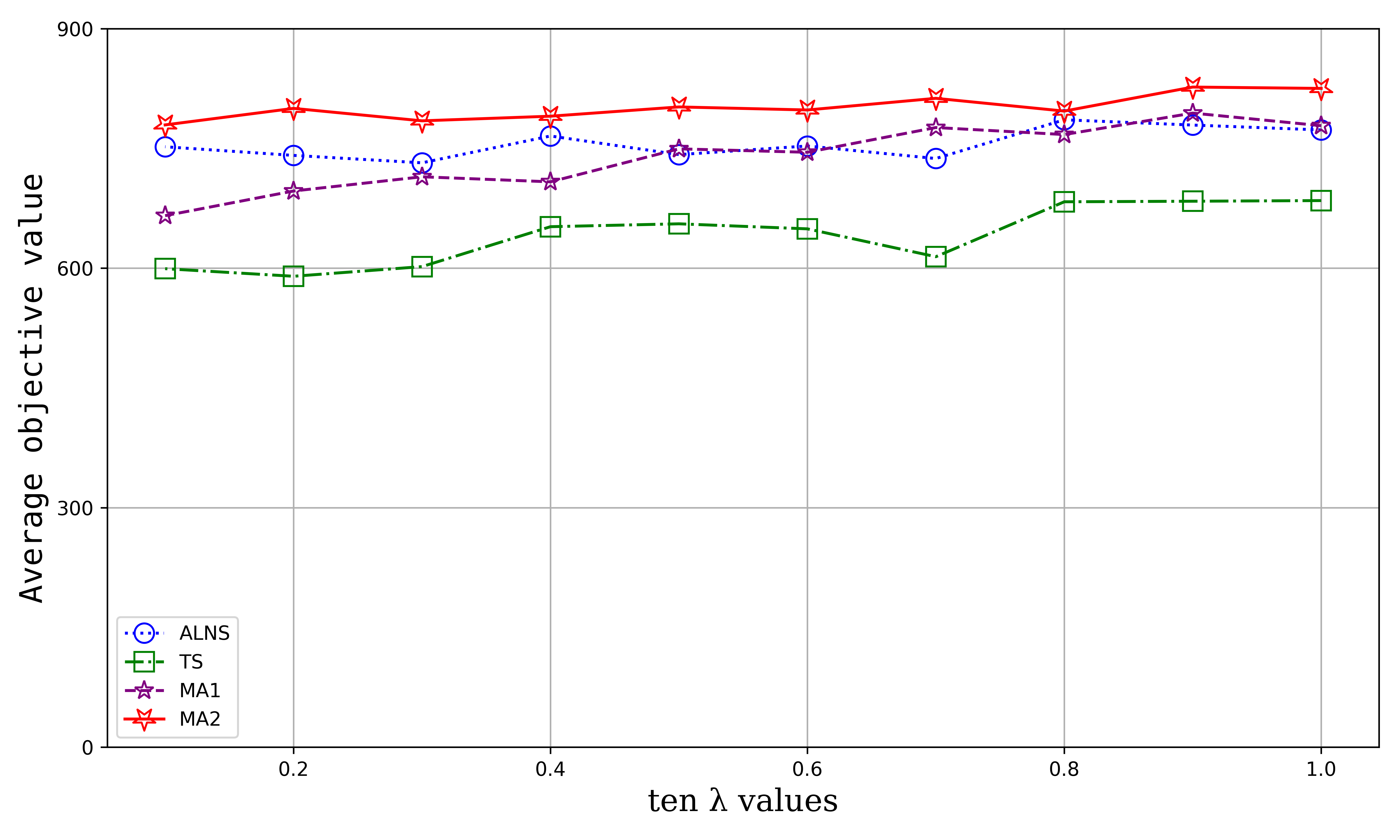}
  \caption{Existing customers=40, New customers=40, $\mu$=0.9}
\end{subfigure}

\begin{subfigure}{\textwidth}
  \centering
  \includegraphics[width=0.65\linewidth]{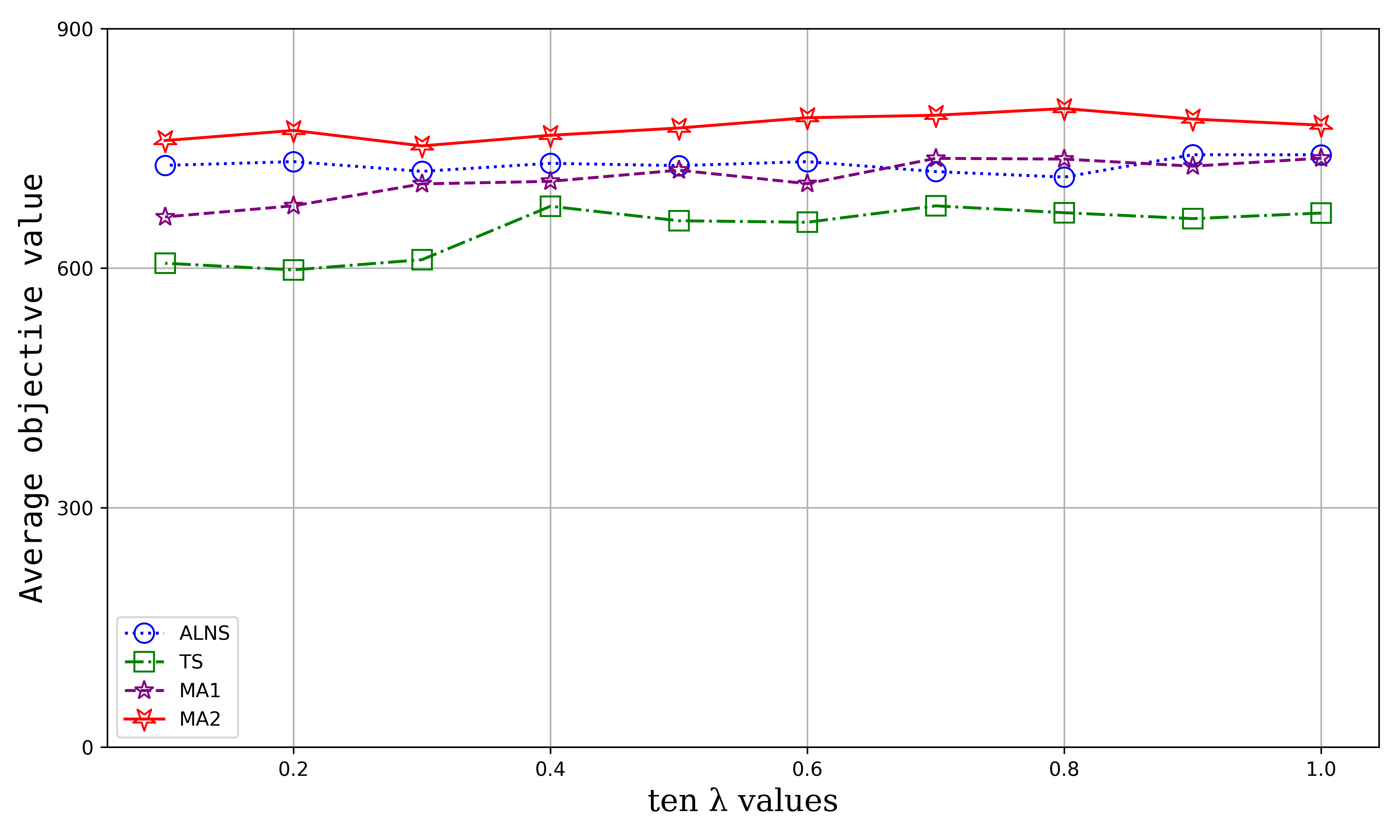}
  \caption{Existing customers=40, New customers=40, $\mu$=0.8}
\end{subfigure}
\caption{Average objective values obtained by MA2, ALNS, TS, and MA1 with 40 existing and 40 new customers}
\label{fig7}
\end{figure}

\begin{figure}[htbp]
\centering
\begin{subfigure}{\textwidth}
  \centering
  \includegraphics[width=0.65\linewidth]{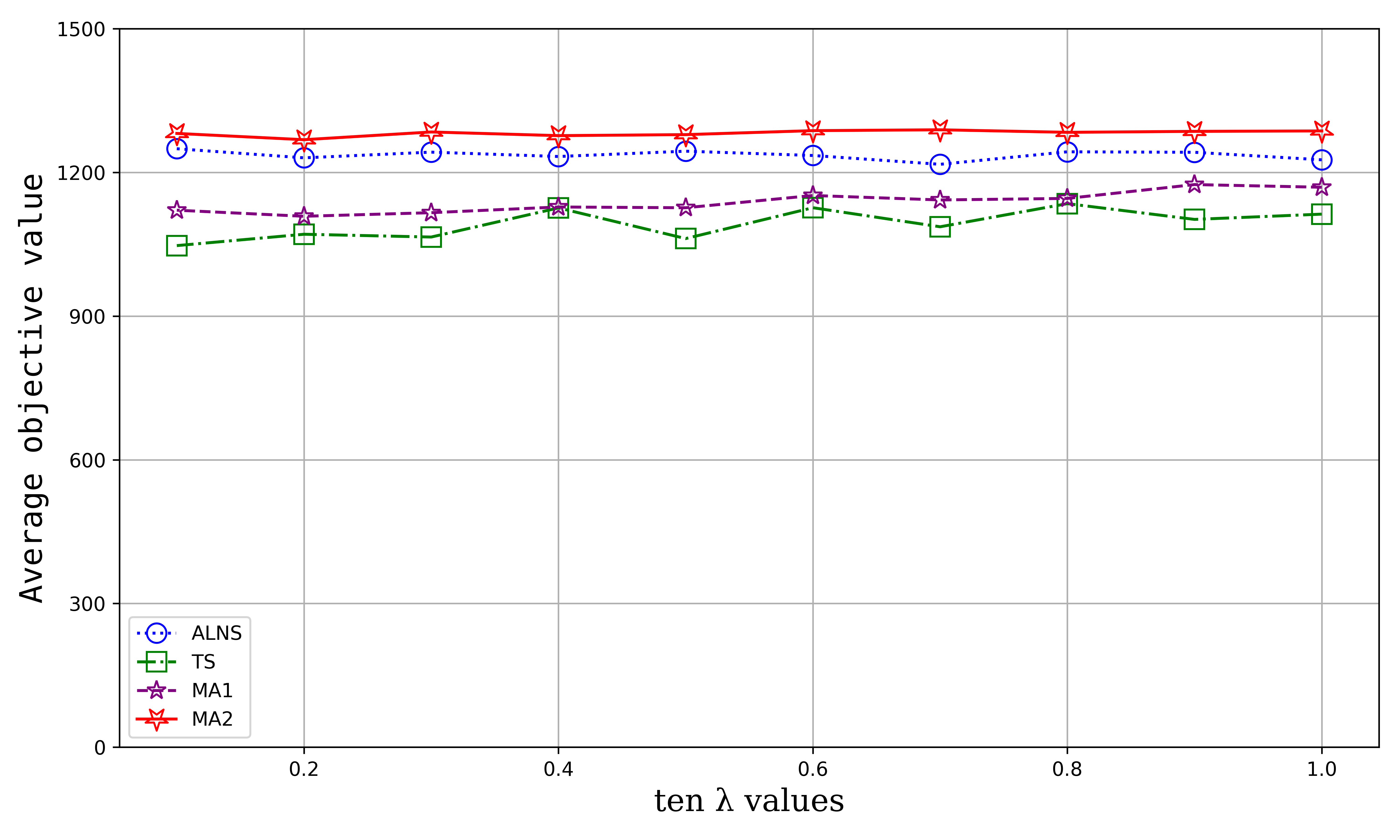}
  \caption{Existing customers=50, New customers=50, $\mu$=1}
\end{subfigure}

\begin{subfigure}{\textwidth}
  \centering
  \includegraphics[width=0.65\linewidth]{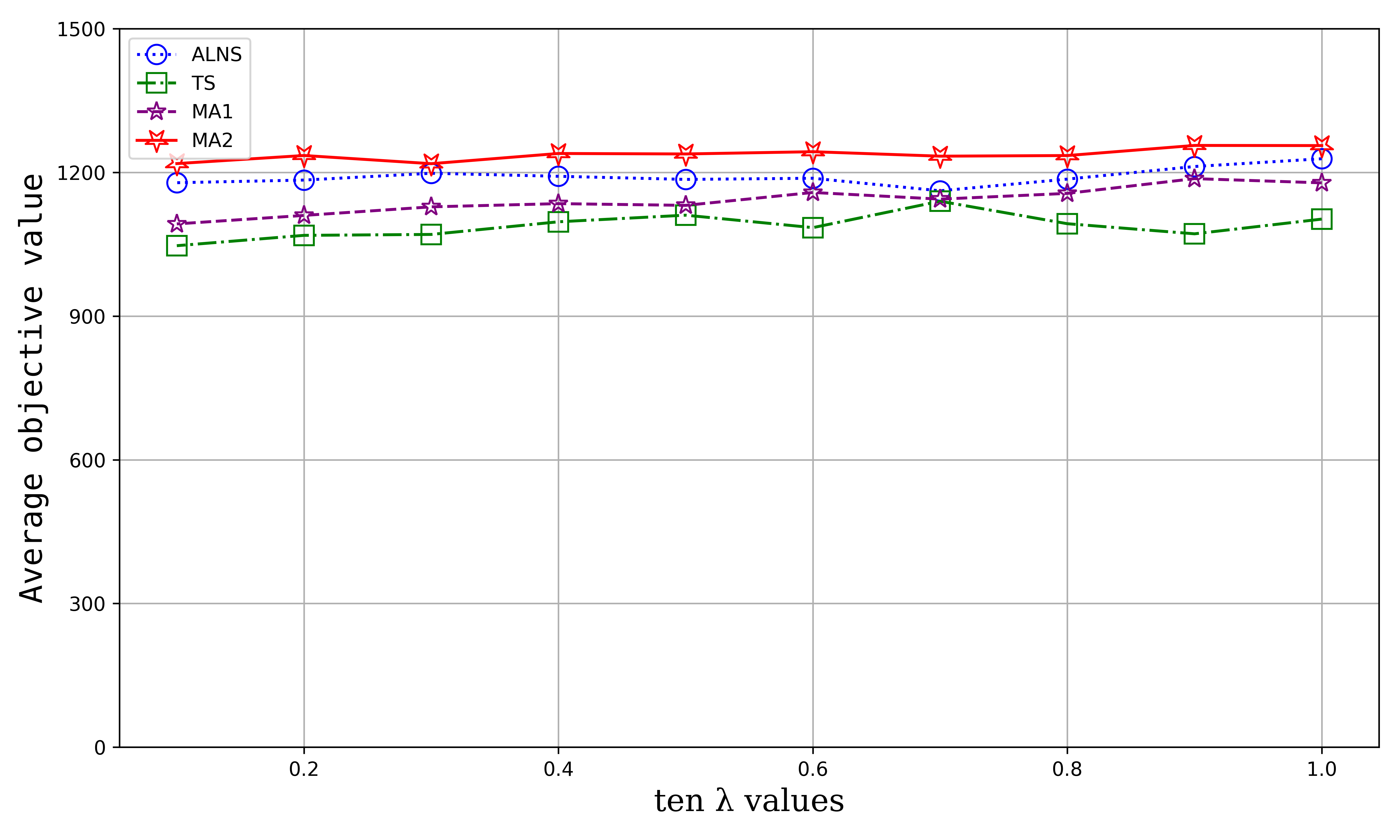}
  \caption{Existing customers=50, New customers=50, $\mu$=0.9}
\end{subfigure}

\begin{subfigure}{\textwidth}
  \centering
  \includegraphics[width=0.65\linewidth]{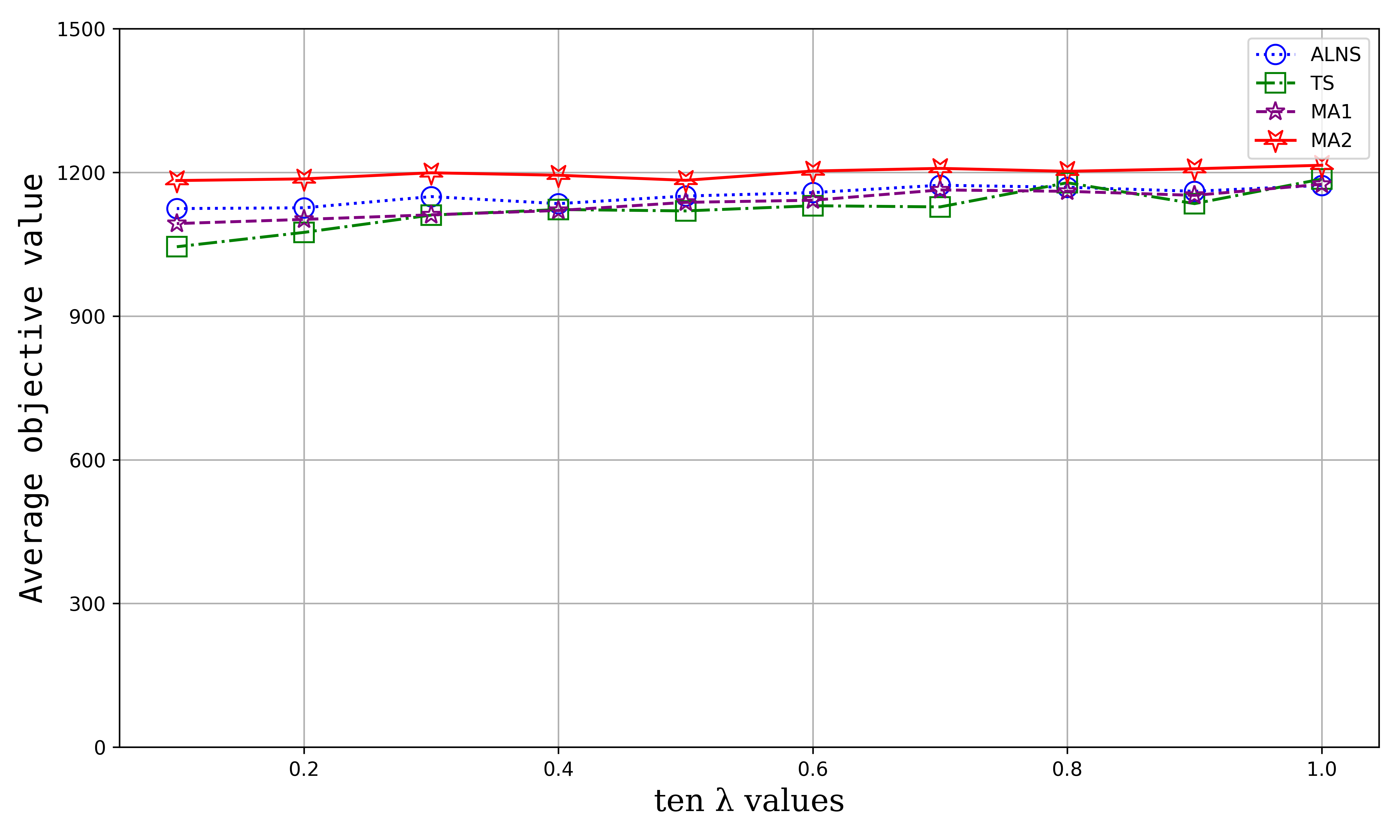}
  \caption{Existing customers=50, New customers=50, $\mu$=0.8}
\end{subfigure}
\caption{Average objective values obtained by MA2, ALNS, TS, and MA1 with 50 existing and 50 new customers}
\label{fig8}
\end{figure}

\begin{landscape}
\begin{figure}[htbp]
\centering
\begin{subfigure}{.31\linewidth}
  \centering
  \includegraphics[width=\linewidth]{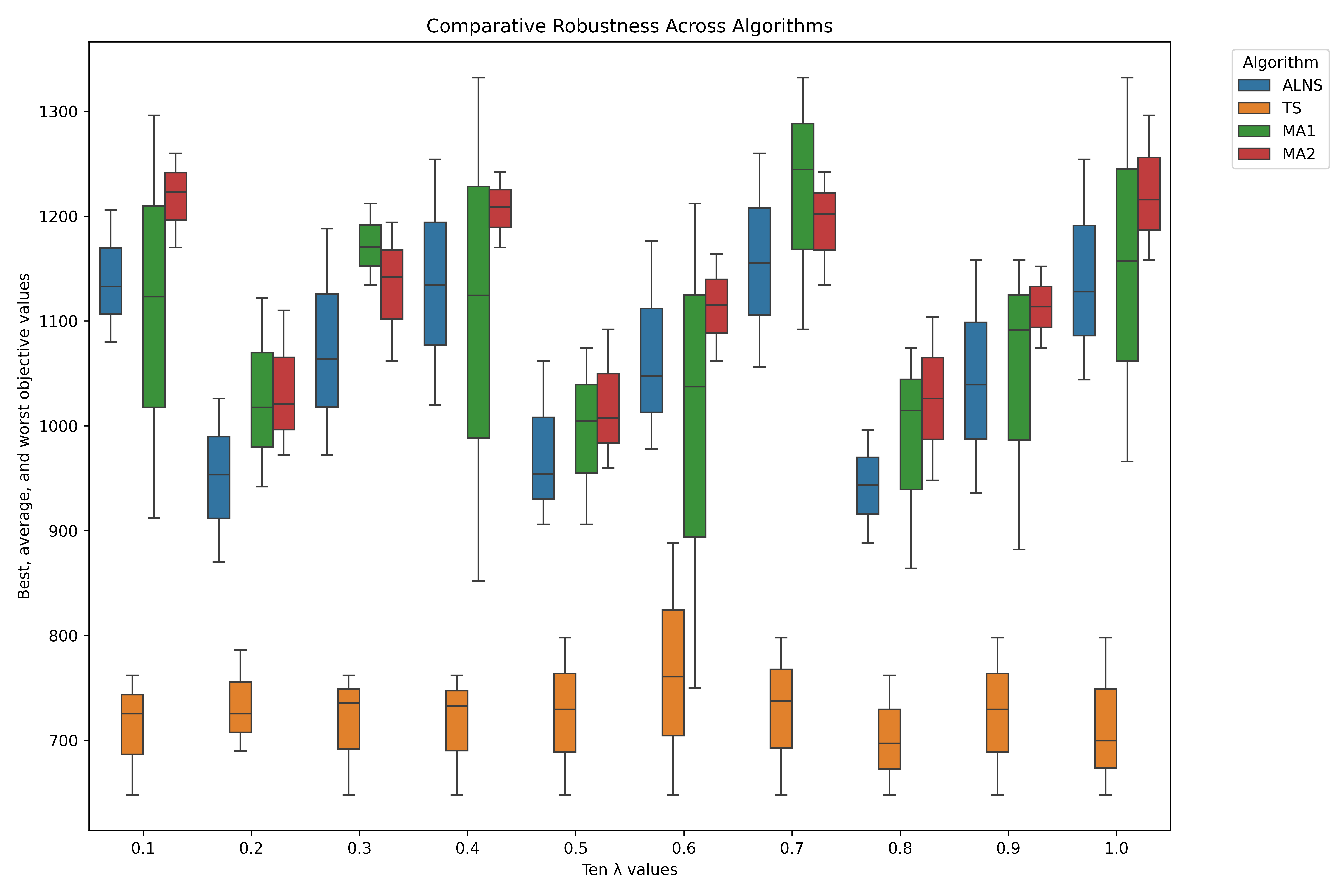}
  \caption{Existing customers=30, New customers=30, $\mu$=1}
\end{subfigure}
\hfill
\begin{subfigure}{.31\linewidth}
  \centering
  \includegraphics[width=\linewidth]{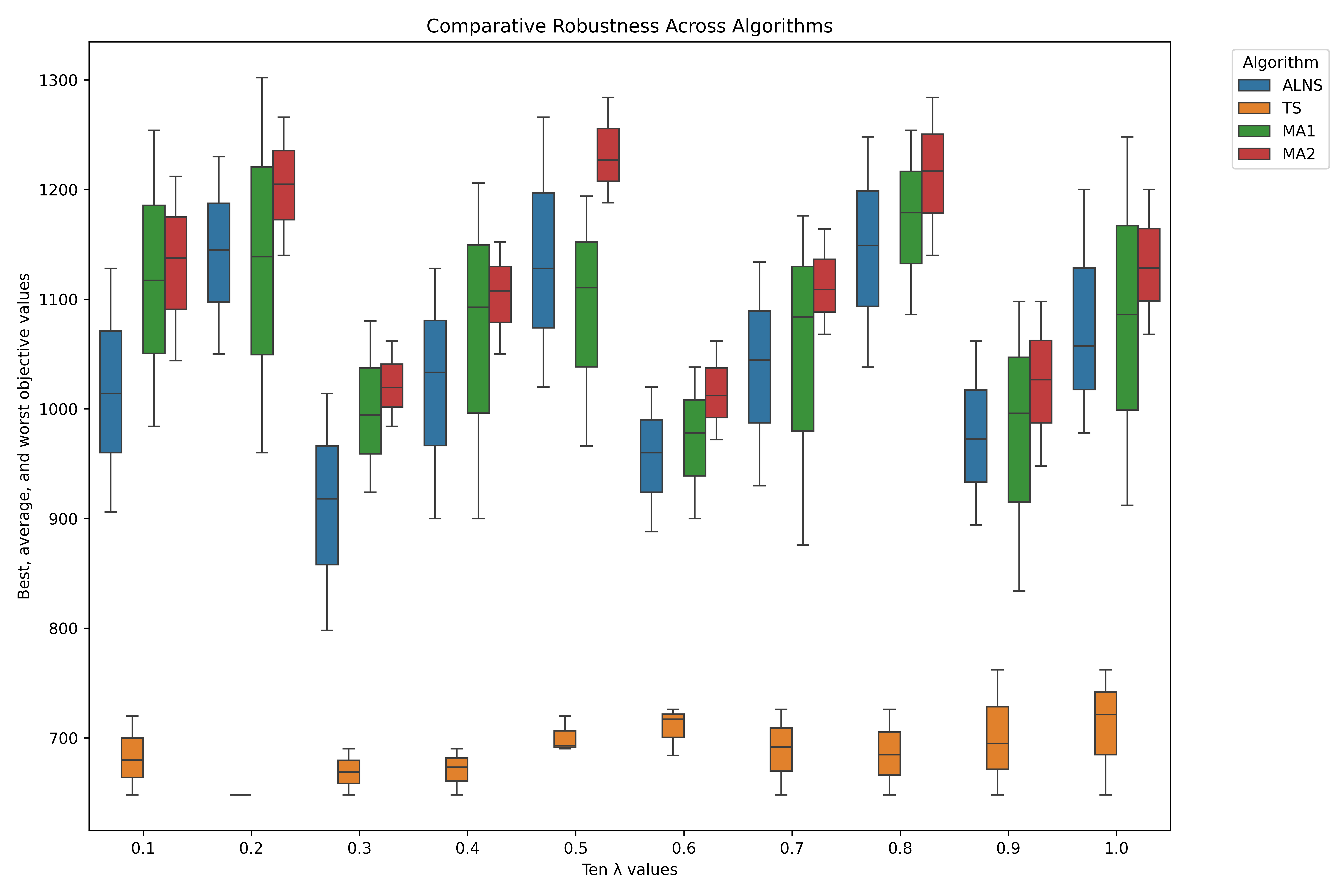}
  \caption{Existing customers=30, New customers=30, $\mu$=0.9}
\end{subfigure}
\hfill
\begin{subfigure}{.31\linewidth}
  \centering
  \includegraphics[width=\linewidth]{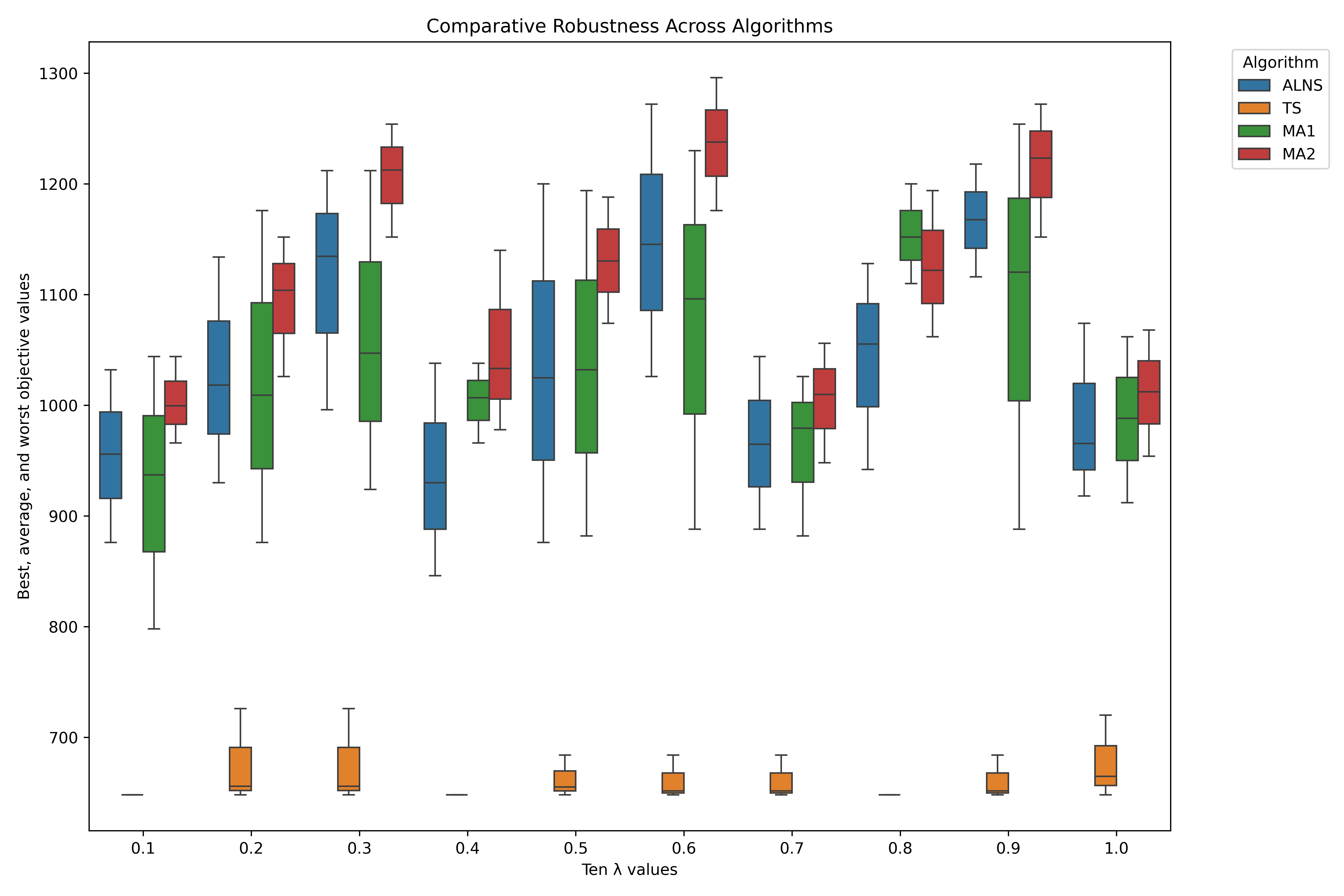}
  \caption{Existing customers=30, New customers=30, $\mu$=0.8}
\end{subfigure}
\newline

\begin{subfigure}{.31\linewidth}
  \centering
  \includegraphics[width=\linewidth]{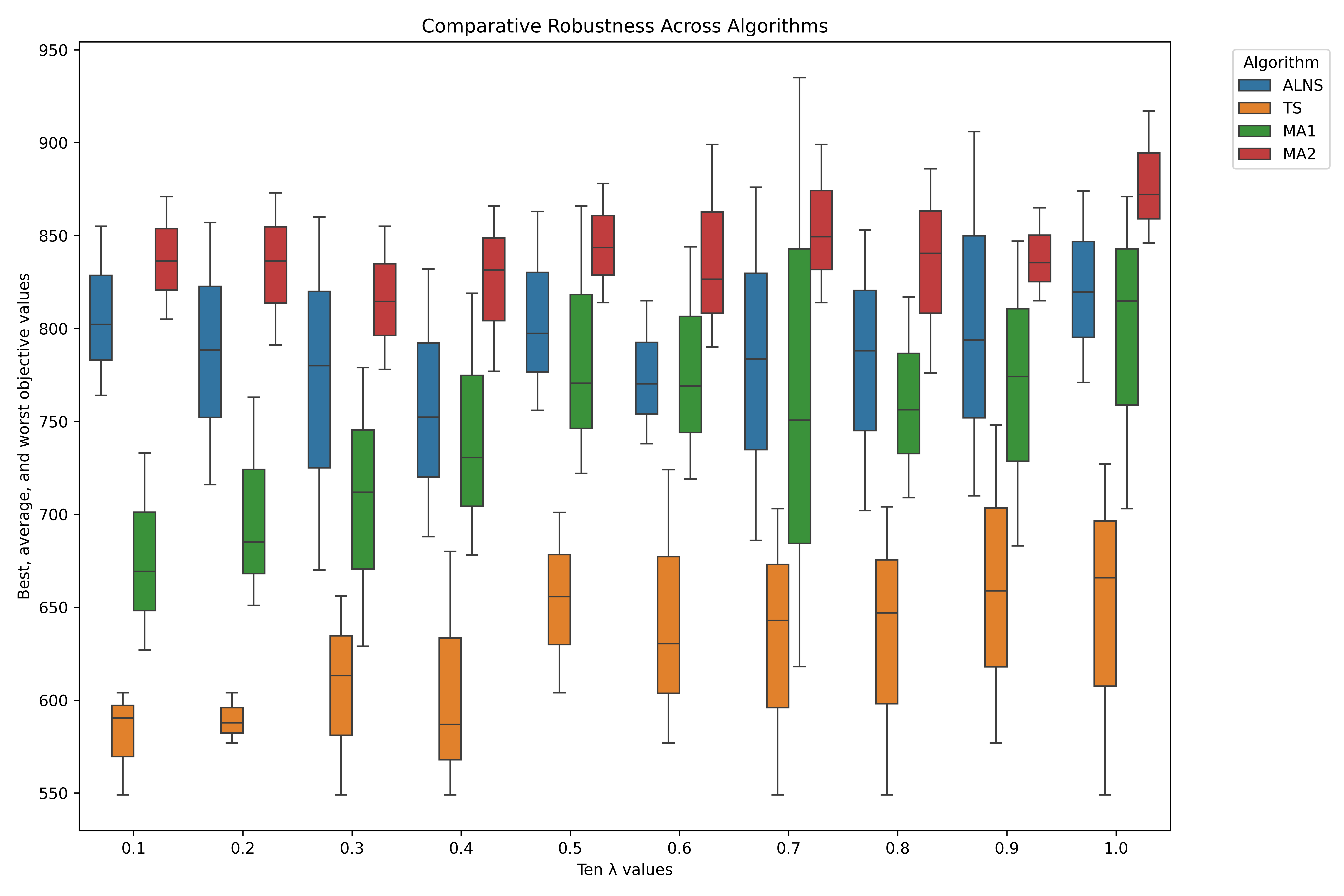}
  \caption{Existing customers=40, New customers=40, $\mu$=1}
\end{subfigure}
\hfill
\begin{subfigure}{.31\linewidth}
  \centering
  \includegraphics[width=\linewidth]{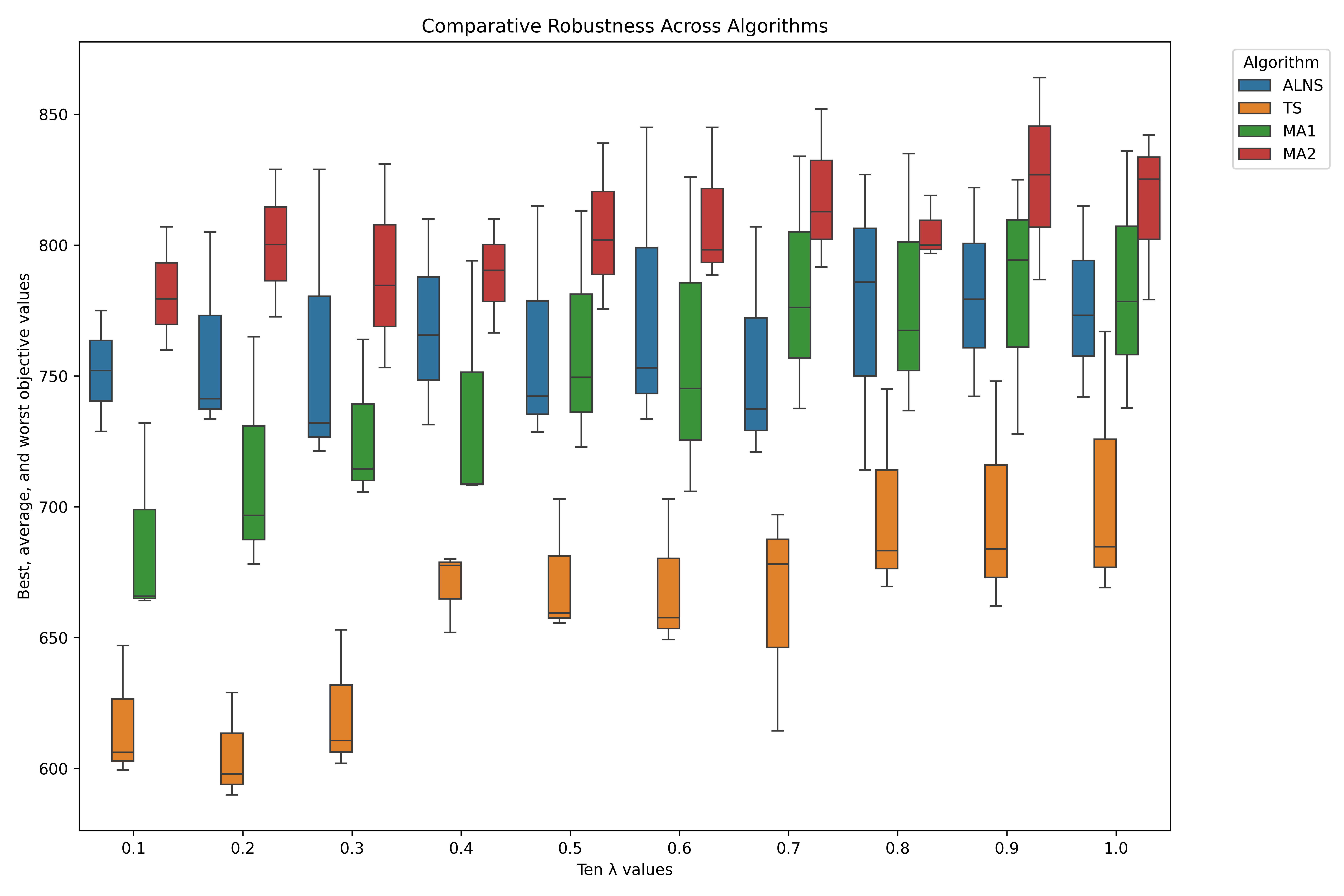}
  \caption{Existing customers=40, New customers=40, $\mu$=0.9}
\end{subfigure}
\hfill
\begin{subfigure}{.31\linewidth}
  \centering
  \includegraphics[width=\linewidth]{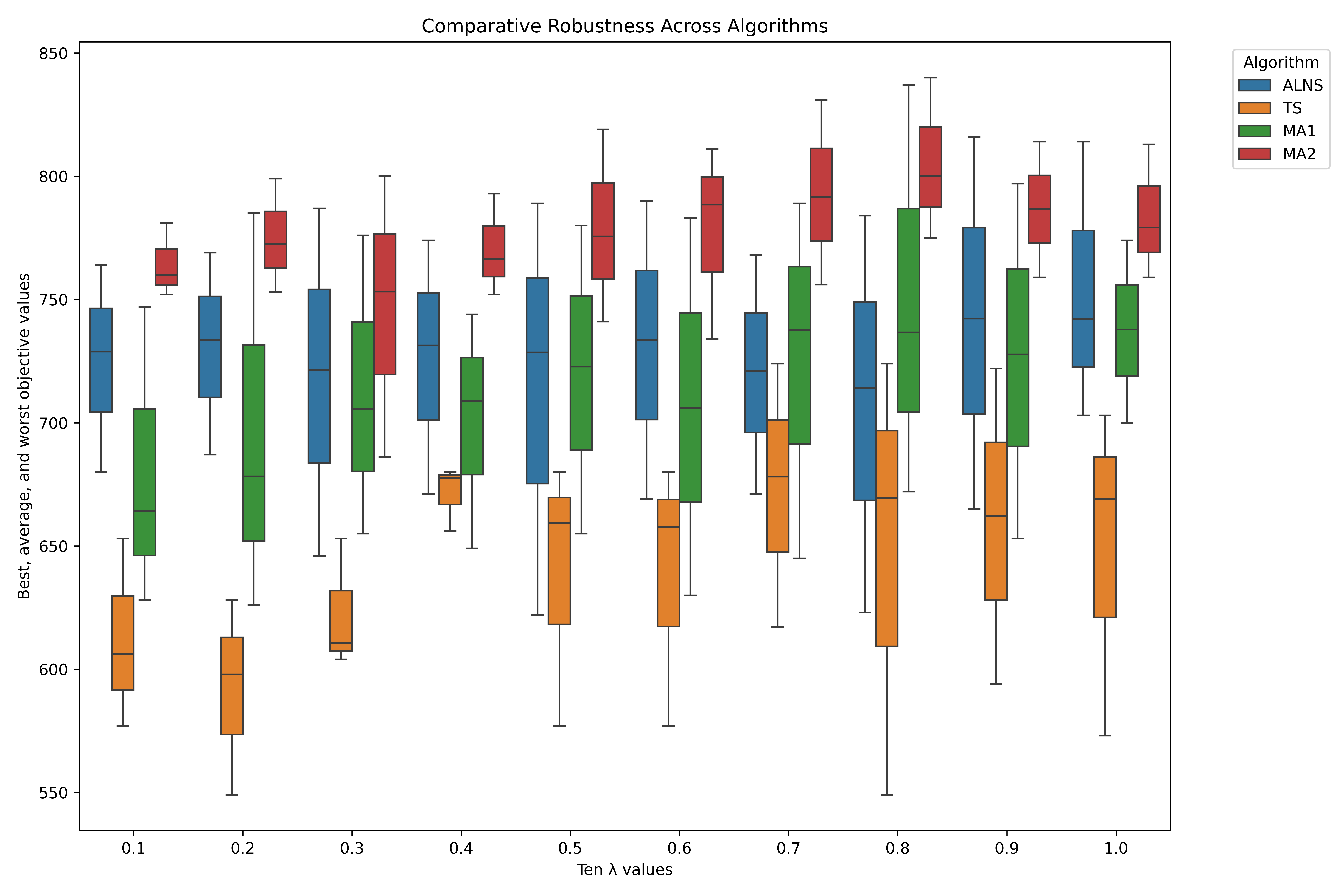}
  \caption{Existing customers=40, New customers=40, $\mu$=0.8}
\end{subfigure}
\newline

\begin{subfigure}{.31\linewidth}
  \centering
  \includegraphics[width=\linewidth]{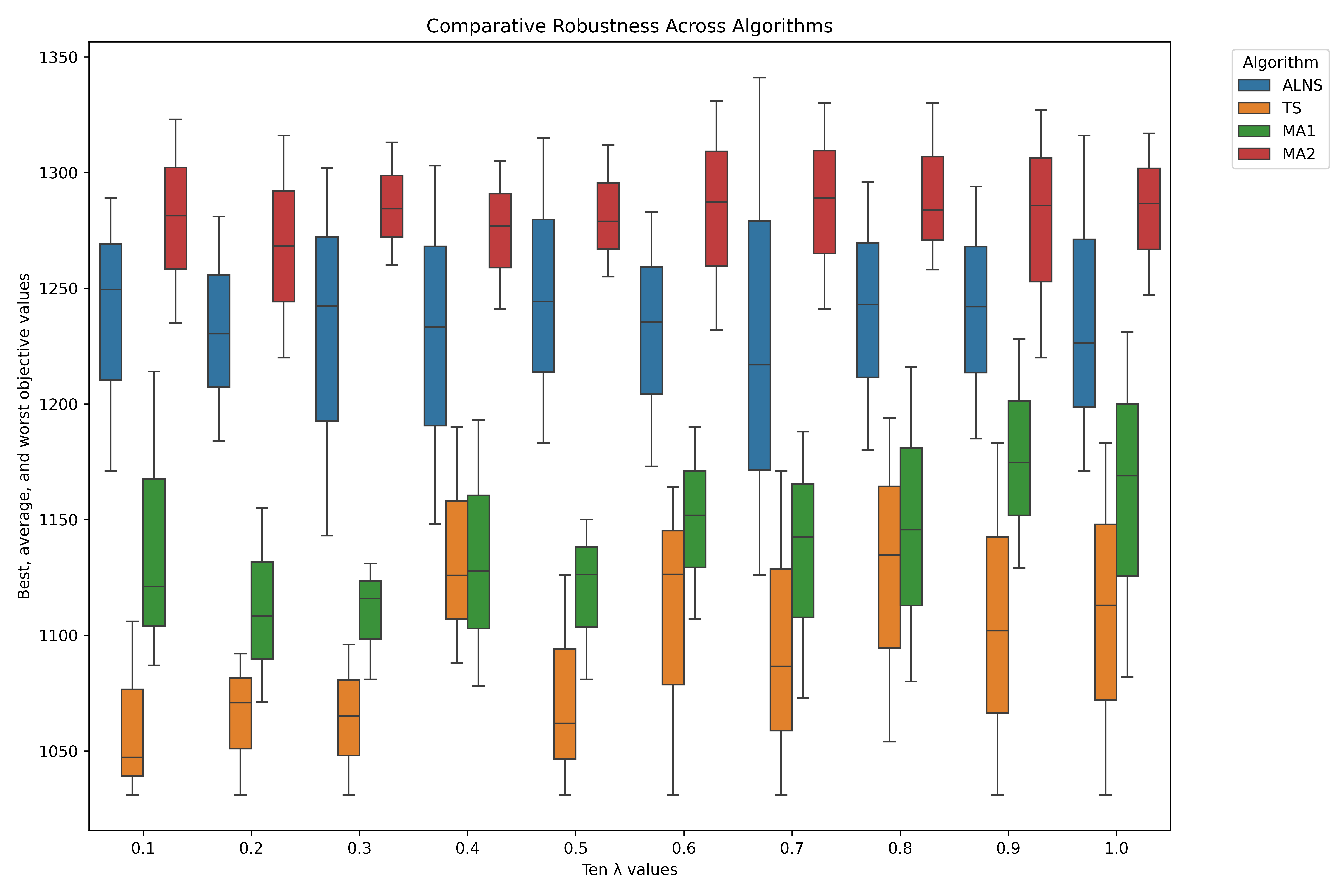}
  \caption{Existing customers=50, New customers=50, $\mu$=1}
\end{subfigure}
\hfill
\begin{subfigure}{.31\linewidth}
  \centering
  \includegraphics[width=\linewidth]{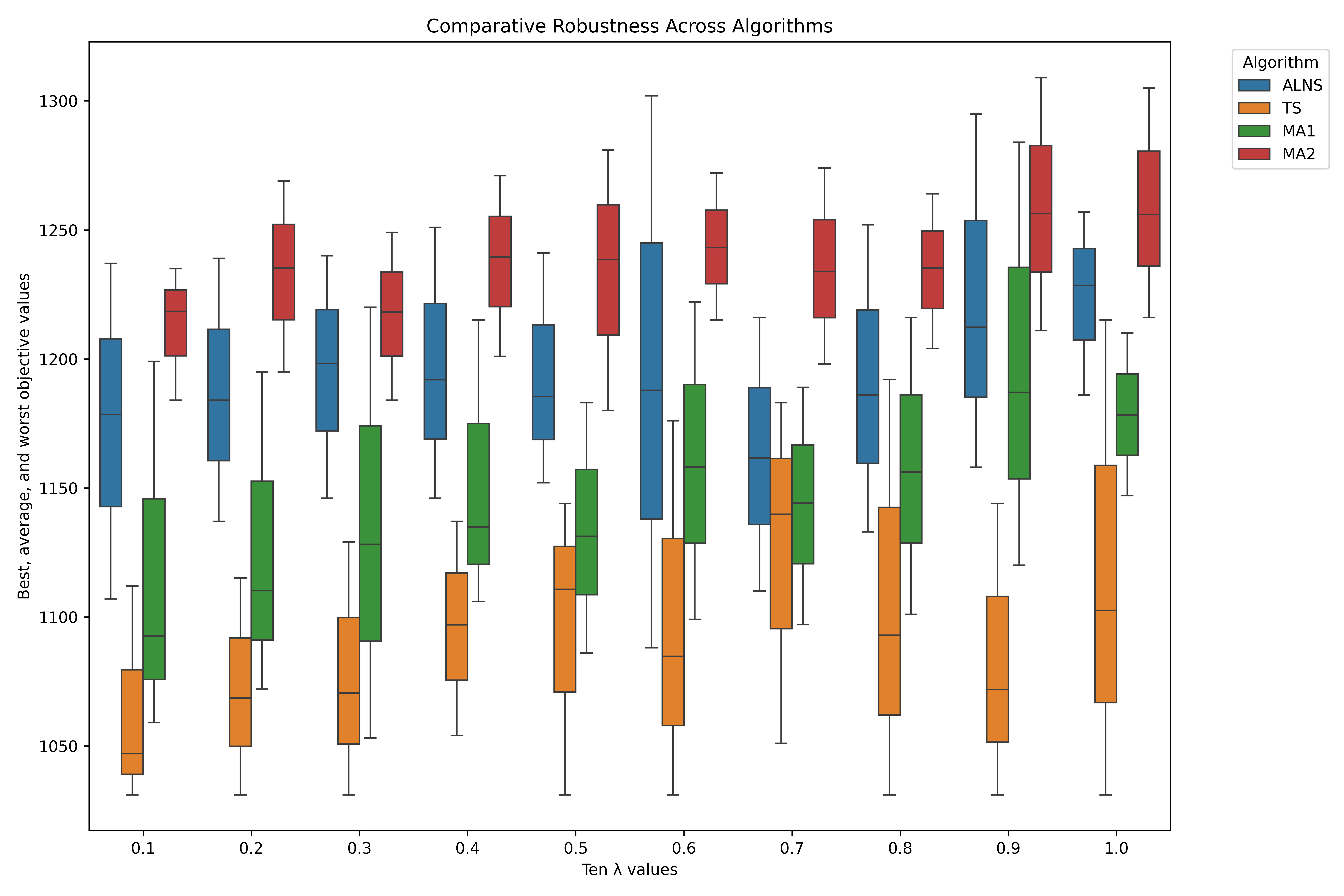}
  \caption{Existing customers=50, New customers=50, $\mu$=0.9}
\end{subfigure}
\hfill
\begin{subfigure}{.31\linewidth}
  \centering
  \includegraphics[width=\linewidth]{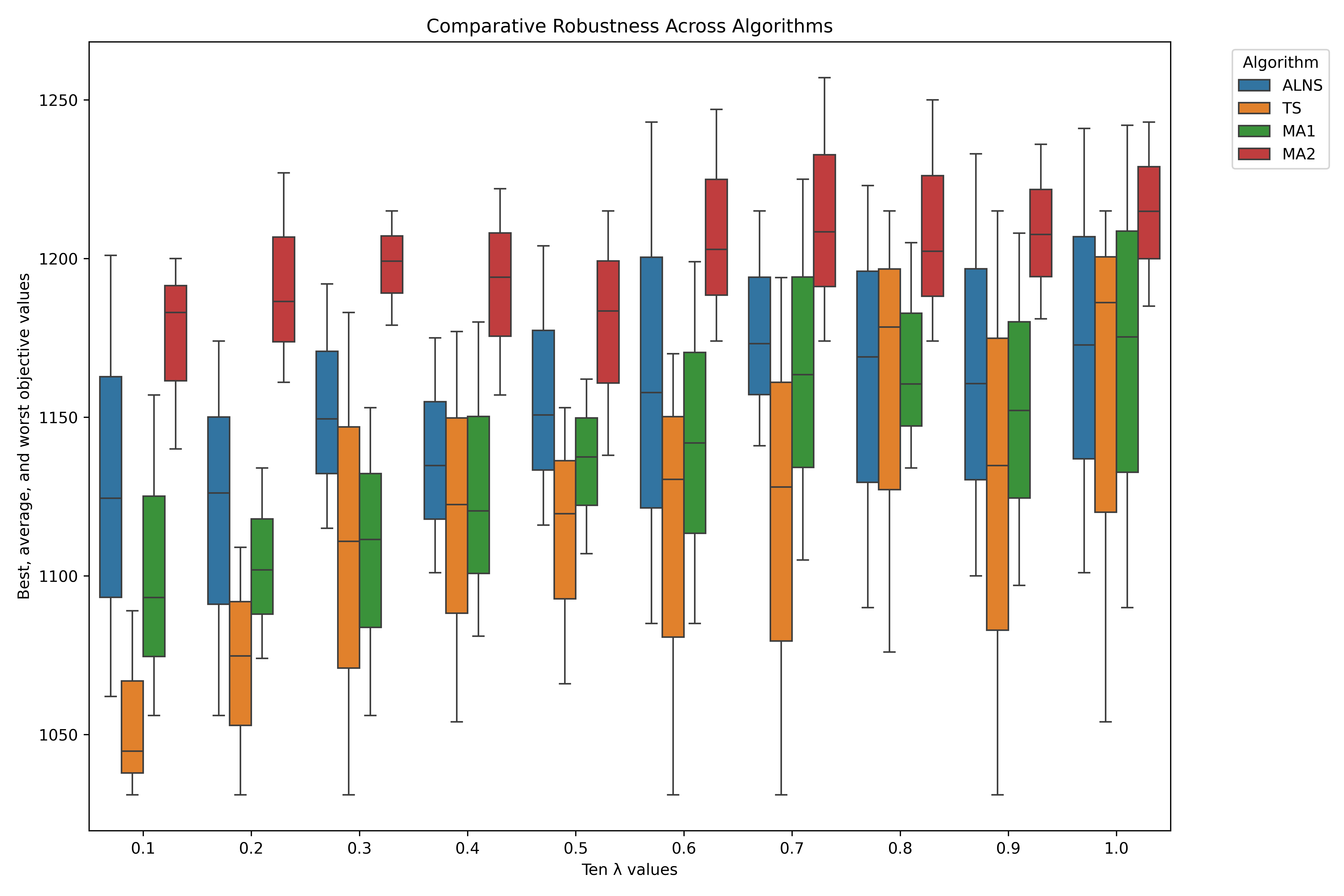}
  \caption{Existing customers=50, New customers=50, $\mu$=0.8}
\end{subfigure}

\caption{Best, average, and worst objective values obtained by MA2, ALNS, TS, and MA1}
\label{fig9}
\end{figure}
\end{landscape}

Detailed comparisons across the four algorithms per instance are presented in the Appendix section. In general, we can firstly observe in Figure \ref{fig6} to Figure \ref{fig8} that MA2 outperforms the others in terms of average performance on 87 of the 90 instances, covering 96.67 of all cases. Next, Figure \ref{fig9} presents box-plot graphs comparing the distributions of the best, average, and worst objective values achieved by the four algorithms. From them, we can conclude that MA2 guarantees the highest level of robustness.

\vspace*{-2cm}

\section{Conclusions and future research}\label{sec_conclusion}
The study successfully developed and evaluated the MA2 algorithm, a hybrid memetic-adaptive neighborhood search optimization approach, for the HHC\&HCRRP-RNC. The algorithm's design specifically addresses the scheduling and routing challenges within the home health care and home care sector, offering a robust solution in the face of uncertainties and dynamic service requests. Computational experiments highlighted MA2's effectiveness in improving operational efficiency, reducing costs, and enhancing service quality compared to traditional methods. By accommodating new customer requests and optimizing caregivers' routes with minimal disruption, MA2 contributes significantly to the resilience and adaptability of home health care services. This research not only advances the application of optimization algorithms in healthcare logistics but also lays the groundwork for future studies to further refine and adapt such algorithms for broader healthcare system challenges.

\bibliographystyle{unsrtnat}
\bibliography{references}  






\end{document}